\documentclass[leqno,a4paper,12pt]{article}
\usepackage{amsfonts,amssymb}
\font\secal=eusm10
\font\ecal=eusm10 at 12pt
\def\esm#1{\hbox{\ecal {#1}}}

\def\Z{\mathbb{Z}}
\def\Q{\mathbb{Q}}

\def\C{\mathbb{C}}

\def\H{\mathbb{H}}
\def\P{\mathbb{P}}

\def\O{{\cal O}}

\def\til#1{\widetilde{#1}}
\def\ovl#1{\overline{#1}}
\def\pf{{\indent\textit{Proof.}\ }}
\def\qed{\hfill$\square$}
\def\Aut{\mathrm{Aut}}

\def\GL{\mathrm{GL}}
\def\PGL{\mathrm{PGL}}
\def\tree{\esm{T}}\def\stree{\hbox{\secal T}}

\def\limitset{\esm{L}}
\def\fixset{\esm{F}}
\def\valuation{\nu}
\def\max{\mathrm{max}}
\def\an{\mathrm{an}}
\def\orb{\mathrm{orb}}
\def\Vert{\mathrm{Vert}}
\def\Edge{\mathrm{Edge}}
\def\Ends{\mathrm{Ends}}
\def\Star{\mathrm{Star}}
\def\longhookrightarrow{\lhook\joinrel\longrightarrow}

\def\bigdownarrow{\vphantom{\bigg|}\Big\downarrow}
\def\underrel#1#2{\mathrel{\mathop{#1}\limits_{#2}}}
\newcounter{para}[section]
\renewcommand{\thepara}{\thesection.\arabic{para}}
\renewcommand{\thesection}{\arabic{section}}
\renewcommand{\paragraph}{\refstepcounter{para}
\indent{\bf\thepara.}}
\newcommand{\sectioning}{\refstepcounter{section}
\indent{\bf \thesection.}}

\newenvironment{condlist}{\vspace{1ex}
     \begin{list}{}{%
          \setlength{\topsep}{0pt}%
          \setlength{\parsep}{0pt}%
          \setlength{\itemsep}{0pt}}%
     }{\end{list}\vspace{1ex}}
\def\cross{P}
\def\msq23{\sqrt{6}}
\def\m2sq3{2\sqrt{3}}
\def\twothree{\heartsuit}
\def\threefour{\diamondsuit}
\title{{\large {\bf $p$-adic Schwarzian triangle groups of Mumford type}}}
\author{Fumiharu Kato}
\date{}
\begin{document}
\maketitle
\begin{center}
\begin{minipage}{10.5cm}
\setlength{\baselineskip}{.85\baselineskip}
\begin{small}
We define a certain $p$-adic analogue of classical Schwarzian
triangle groups related to Mumford's uniformization of analytic
curves and give a complete classification of it. 
\end{small}
\end{minipage}
\end{center}

\vspace{1ex}
\sectioning\label{section-introduction}\ 
{\bf Introduction}

\vspace{1ex}
\paragraph\label{para-orbifolds}\ 
{\sl Uniformization of orbifolds and triangle groups.}\ 
The rich geometric structure of uniformized analytic varieties over
non-archimedean fields has been studied by many authors, and already has
a long history. 
Mumford \cite{Mum72} showed that an analytic curve $\mathcal{X}$ defined
over a non-archimedean field $K$ with a split multiplicative analytic
reduction can be uniformized as $\mathcal{X}\cong\Gamma\backslash
(\P^{1,\an}_K-\limitset_{\Gamma})$, where 
$\Gamma$ is a finitely generated free discrete subgroup of $\PGL(2,K)$
and $\limitset_{\Gamma}$ is the set of limit points of $\Gamma$.
An equally important example is the uniformization of a curve which is
an \'etale covering of a Mumford curve, studied by van der Put
\cite{vdP83}.  
These are the most practical and reasonable analogues of uniformizations
of complex analytic curves. 

Historically, however, the theory of uniformization in complex analysis
arose from interplay between geometric and function-theoretic viewpoints. 
This is apparent if one considers the orbifold uniformization of 
$\P^{1,\an}_{\C}$ with finitely many points marked by positive
integers ($=$ the ramification degrees). 
The link between geometry and function theory stems from the behavior of 
the multivalued function $z$, inverse to the uniformization map, which
is written as a ratio $z=u_1/u_2$ of two linearly independent solutions
of a Fuchsian differential equation with rational exponents (cf.\
\cite{Yos87}).  
For example, if the orbifold is $\P^{1,\an}_{\C}$ with
precisely the three points $0,1,\infty$ marked by
$e_0,e_1,e_{\infty}\in\Z_{>0}$, then the corresponding
differential equation is the Gaussian hypergeometric equation
$$
x(1-x)\frac{d^2u}{dx^2}+\{c-(a+b+1)x\}\frac{du}{dx}-abu=0
\leqno{\indent\textrm{(\ref{para-orbifolds}.1)}}
$$
(where $a,b,c\in\Q$ with $|1-c|=1/e_0$, $|c-a-b|=1/e_1$, and 
$|a-b|=1/e_{\infty}$).
Depending on whether $1/e_0+1/e_1+1/e_{\infty}>1$, $=1$ or $<1$, the
image of $z=u_1/u_2$ is isomorphic to $\P^{1,\an}_{\C}$, $\C$ or
$\H$, respectively, and $z$ maps the upper-half plane onto the interior of a
triangle region with the angles $\pi/e_0$, $\pi/e_1$ and $\pi/e_{\infty}$.  
The corresponding orbifold fundamental group has a representation into
$\PGL(2,\C)$ with discrete image, the so-called {\it Schwarzian
triangle group} $\Delta(e_0,e_1,e_{\infty})$, whose action on the
universal covering space is visible in terms of complex reflections (cf.\ 
\cite[Chap.\ II]{Mag74}). 

The problems which arises from carrying out such a program in the
$p$-adic situation are mainly topological. It is perhaps appropriate here
to remind the reader of the fact that in rigid analysis an \'etale
covering map is not necessarily a topological covering map ($=$ the
locally topologically trivial map). Consequently, in contrast to complex 
analysis, we have many simply connected domains, even one dimensional;
for instance, the complement of finitely many points in
$\P^{1,\an}_K$ is simply connected. In particular, a reasonable
definition of the orbifold fundamental groups is highly non-trivial.

In \cite{And98}, Y.\ Andr\'e studied the rather special class of \'etale
covering maps which are composites of topological coverings followed
by finite \'etale (not necessarily topological) coverings. He observed
that the covering maps of this kind give rise to a reasonable concept
of orbifold fundamental groups (denoted by
$\pi^{\orb}_1$ in \cite{And98}) in the $p$-adic situation, and
discussed the relation with differential equations; one of his result
phrases it as follows \cite[\S 6]{And98}: Consider the orbifold
$\mathcal{X}$ (cf.\ \cite[5.1]{And98} for the precise definition) which
is supported on 
$\P^{1,\an}_{\C_p}$ with $n$-points $\zeta_i$ marked by positive
integers $e_i$ ($1\leq i\leq n$). Then:  

\vspace{1ex}\noindent
{\sl There exists a canonical fully faithful functor of categories}
$$
\left\{
\begin{minipage}{5cm}
\setlength{\baselineskip}{.85\baselineskip}
\begin{small}
{\slshape Continuous
representations $\rho$ of $\pi^{\orb}_1(\mathcal{X},\ovl{x})$ into 
$\GL(r,\C_p)$ with discrete coimage}
\end{small}
\end{minipage}
\right\}
\longrightarrow
\left\{
\begin{minipage}{5cm}
\setlength{\baselineskip}{.85\baselineskip}
\begin{small}
{\slshape
Algebraic regular connections on $\P^1_{\C_p}-\{\zeta_i\,|\,1\leq i\leq
n\}$ of rank $r$ such that the local monodromy at each $\zeta_i$ is of 
finite order dividing $e_i$}
\end{small}
\end{minipage}
\right\}.
\leqno{\indent\textrm{(\ref{para-orbifolds}.2)}}
$$
{\slshape Moreover, the essential image of this functor consists of the
connections $\nabla$ enjoying the following condition (Global
Monodromy Condition):}

\vspace{1ex}
(\ref{para-orbifolds}.3)\ 
{\sl There exists a connected rigid analytic curve
           $\mathcal{S}$ and a finite morphism
           $\varphi\colon\mathcal{S}\rightarrow\mathcal{X}$ ramified
           above precisely the points $\zeta_i$ with ramification
           indices dividing $e_i$, such that the connection
           $\varphi^{\ast}\nabla$ on $\mathcal{S}$ admits a full set of
           multivalued analytic solutions on $\mathcal{S}$ (and,
           moreover, on the Berkovich space associated to $\mathcal{S}$).}

\vspace{1ex}
In particular, if $n=3$, $r=2$ and if the image $\Gamma$ of $\rho$ in
$\PGL(2,\C_p)$ is discrete, then $\Gamma$ can be regarded as a $p$-adic 
analogue of the Schwarzian triangle group $\Delta(e_0,e_1,e_{\infty})$.
If so, $\Gamma$ gives the ``projective monodromy'' for the connection
$\nabla$ defined by the functor (\ref{para-orbifolds}.2), which is 
nothing but the one associated to 
the Gaussian hypergeometric equation (\ref{para-orbifolds}.1).
In \cite[\S 9]{And98}, Andr\'e discussed such groups, called {\it
$p$-adic triangle groups}, and gave a complete list of the so-called
{\it arithmetic} $p$-adic triangle groups, which are constructed through the
Cherednik-Drinfeld theory of uniformization of Shimura curves, starting
from Takeuchi's list of arithmetic triangle groups.  
Notably, he deduced that {\sl there exists no arithmetic $p$-adic
triangle groups for $p>5$}. The uniformizations of the orbifolds
$\mathcal{X}$ corresponding to these groups are given by the Drinfeld
upper-half plane or its \'etale coverings.

\vspace{1ex}
\paragraph\label{para-results}\ 
{\sl Results of this paper.}\ 
In this paper, we will discuss (not necessarily arithmetic) $p$-adic
triangle groups $\Gamma$ as above, especially in the case that
the corresponding uniformization is given by the space 
$\P^{1,\an}_{\C_p}-\limitset_{\Gamma}$.
We call such a group $\Gamma$ a $p$-adic triangle group {\it of Mumford
type}; we can define it in simpler terms (without involving
$\pi^{\orb}_1$) as follows: 

\vspace{1ex}
{\bf Definition.}\ A finitely generated discrete subgroup $\Gamma$ of
$\PGL(2,\C_p)$ is said to be a {\it $p$-adic (Schwarzian) triangle
group of Mumford type} if
$\Gamma\backslash(\P^{1,\an}_{\C_p}-\limitset_{\Gamma})\cong
\P^{1,\an}_{\C_p}$ and the uniformization map
$$
\varpi_{\Gamma}\colon(\P^{1,\an}_{\C_p}-\limitset_{\Gamma})
\longrightarrow\P^{1,\an}_{\C_p}
$$
is ramified above precisely three points.

\vspace{1ex}
If $\Gamma$ is a finite subgroup, then $\limitset_{\Gamma}=\emptyset$,
and the map $\varpi_{\Gamma}$ is the analytification of the algebraic
quotient map $\P^1_{\C_p}\rightarrow\P^1_{\C_p}/\Gamma$. 
Hence the spherical (i.e., $1/e_0+1/e_1+1/e_{\infty}>1$) $p$-adic triangle
groups of Mumford type amount to the same as the classical ones.
Our main theorem gives the complete classification of the other $p$-adic
triangle groups of Mumford type:

\vspace{1ex}
{\bf Theorem.}\ {\slshape {\rm (1)} A non-spherical $p$-adic triangle
group $\Gamma$ of Mumford type of index $(e_0,e_1,e_{\infty})$ exists if
and only if $p$ and the unordered triple $(e_0,e_1,e_{\infty})$ occur in 
the left two columns of Table 1, where, in the last row, the integers
$l,m,n$ obey the following condition;

\vspace{1ex}
{\rm ($\ast$)} $lmn\geq 4$, and at least two of $l,m,n$ are odd.

\vspace{1ex}\noindent
In each of the cases in Table 1, the conjugacy class of
$\Gamma$ in $\PGL(2,\C_p)$ is uniquely determined. 
In particular, $p$-adic triangle groups of Mumford type do not exist if
$p>5$. }

\begin{small}
\begin{center}
{\sc Table 1}: List of $p$-adic triangle groups of Mumford type

\vspace{2ex}
\begin{tabular}{|c|c|cc|}
\hline
&$p$&\multicolumn{2}{|c|}{Index}\\
\hline\hline
I&$5$&$(3,3,5l)$&$(l\geq 1)$\\
\hline
II&$5$&$(2,3,5l)$&$(l\geq 2)$\\
\hline\hline
III&$3$&$(5,5,3l)$&$(l\geq 1)$\\
\hline
IV&$3$&$(2,5,3l)$&$(l\geq 2)$\\
\hline
V&$3$&$(4,4,3l)$&$(l\geq 1)$\\
\hline
VI&$3$&$(2,3l,3l)$&$(l\geq 2)$\\
\hline
VII&$3$&$(2,4,3l)$&$(l\geq 2)$\\
\hline
VIII&$3$&$(4,5,3l)$&$(l\geq 1)$\\
\hline\hline
IX&$2$&$(5,5,2l)$&$(l\geq 1\ \textrm{and}\ l\textrm{:\ odd})$\\
\hline
X&$2$&$(5,5,4l)$&$(l\geq 1)$\\
\hline
XI&$2$&$(5,2l,4m)$&$(lm\geq 1\ \textrm{and}\ l\textrm{:\ odd})$\\
\hline
XII&$2$&$(3,3,2l)$&$(l\geq 3\ \textrm{and}\ l\textrm{:\ odd})$\\
\hline
XIII&$2$&$(3,3,4l)$&$(l\geq 1)$\\
\hline
XIV&$2$&$(3,2l,4m)$&$(lm\geq 2\ \textrm{and}\ l\textrm{:\ odd})$\\
\hline
XV&$2$&$(3,5,2l)$&$(l\geq 3\ \textrm{and}\ l\textrm{:\ odd})$\\
\hline
XVI&$2$&$(3,5,4l)$&$(l\geq 1)$\\
\hline
XVII&$2$&$(2m+1,2l,2l)$&$(m\geq 1,l\geq 2)$\\
\hline
XVIII&$2$&$(2l,2m,2n)$&$(l,m,n\ \textrm{satisfy}\ \textrm{($\ast$)})$\\
\hline
\end{tabular}
\end{center}
\end{small}

{\rm (2)} {\sl In each of the cases in Table 1, $\Gamma$ is
isomorphic to the abstract group given as follows (here, $Z_n$ denotes
the cyclic group of order $n$, $D_n$ the dihedral group of 
degree $n$ $(\cong Z_n\rtimes Z_2)$, etc. The symbol $\ast$ means the
amalgam product, that is, the push-forward in the category of groups.):}

$$
\begin{array}{cll}
\textrm{I}&\textrm{:}&A_5\ast_{D_5}D_{5l}\ast_{D_5}A_5,\\
\\
\textrm{II}&\textrm{:}&A_5\ast_{D_5}D_{5l},\\
\\
\textrm{III}&\textrm{:}&A_5\ast_{D_3}D_{3l}\ast_{D_3}A_5,\\
\\
\textrm{IV}&\textrm{:}&A_5\ast_{D_3}D_{3l},\\
\\
\textrm{V}&\textrm{:}&S_4\ast_{D_3}D_{3l}\ast_{D_3}S_4,\\
\\
\textrm{VI}&\textrm{:}&A_4\ast_{Z_3}Z_{3l},\\
\\
\textrm{VII}&\textrm{:}&S_4\ast_{D_3}D_{3l},\\
\\
\textrm{VIII}&\textrm{:}&S_4\ast_{D_3}D_{3l}\ast_{D_3}A_5,\\
\\
&&\qquad\quad D_{2l}\\
\textrm{IX}&\textrm{:}&\qquad\quad\ast_{D_2}\\
&&A_5\ast_{A_4}A_4\ast_{A_4}A_5,\\
\end{array}
\qquad\qquad
\begin{array}{cll}
\textrm{X}&\textrm{:}&A_5\ast_{A_4}S_4\ast_{D_4}D_{4l}\ast_{D_4}S_4
\ast_{A_4}A_5,\\
\\
&&\qquad\qquad A_5\\
\textrm{XI}&\textrm{:}&\qquad\qquad \ast_{A_4}\\
&&D_{4m}\ast_{D_4}S_4\ast_{D_2}D_{2l},\\
\\
\textrm{XII}&\textrm{:}&A_4\ast_{D_2}D_{2l},\\
\\
\textrm{XIII}&\textrm{:}&S_4\ast_{D_4}D_{4l}\ast_{D_4}S_4,\\
\\
\textrm{XIV}&\textrm{:}&D_{4m}\ast_{D_4}S_4\ast_{D_2}D_{2l},\\
\\
\textrm{XV}&\textrm{:}&A_5\ast_{D_2}D_{2l},\\
\\
\textrm{XVI}&\textrm{:}&A_5\ast_{A_4}S_4\ast_{D_4}D_{4l}\ast_{D_4}S_4,\\
\\
\textrm{XVII}&\textrm{:}&D_{2m+1}\ast_{Z_2}Z_{2l},\\
\\
\textrm{XVIII}&\textrm{:}&D_{2l}\ast_{D_2}D_{2m}\ast_{D_2}D_{2n}.
\end{array}
$$

\vspace{1ex}\noindent
{\sl Here, in I, III, V, VIII, the two dihedral
groups with the same order are chosen to be equal if $l$ is odd, and are 
not equal and are conjugate with each other by an involution in the
dihedral group between denoted between them if $m$ is even.
In X, XIII, XVI, and XVIII, the two dihedral
groups with the same order are chosen to be equal, and in XI and XIV, the
subgroups $D_4$ and $D_2$ in $S_4$ are chosen such that $D_2\subset
D_4$.
In XI, the subgroup $A_4$ in $S_4$ is chosen so that its intersection
with $D_2$ is trivial.}

\vspace{1ex}
{\bf Remarks.}\ (1)\ 
The theorem, in particular, proves Yves Andr\'e's
conjecture that there are infinitely many non-arithmetic $p$-adic
triangle groups (cf.\ \cite{And98}). 

(2) In each of the cases in the table, at least one of
the numbers $e_i$ $(i=0,1,\infty)$ is divisible by the residue
characteristic $p$.  

(3) The theorem shows that there are no Euclidean (i.e.,
$1/e_0+1/e_1+1/e_{\infty}=1$) $p$-adic triangle groups of Mumford
type. The reason for this is that the elliptic curve which covers an
Euclidean orbifold always has a complex multiplcation and never be a
Tate curve. 

(4) The theorem and the well-known fact on automorphisms of Mumford
curves \cite[VII.\S 1]{GvP80} imply that for $p>5$ and for a Mumford
curve $X$, if $\Aut(X)\backslash X\cong\P^1_K$, then the quotient map 
$X\rightarrow\P^1_K$ ramifies above at least 4 points.
Applying the classical Hurwitz formula, we therefore see that
$|\Aut(X)|\leq 12(g-1)$ (where $g$ is the genus of $X$), which 
partly recovers the Herrlich's result \cite{Her80b}.

(4) Our list of $p$-adic triangle groups has a non-empty intersection
with Andr\'e's list of $p$-adic arithmetic triangle groups, but does 
not include it, since Andr\'e's $p$-adic arithmetic triangle groups are
not of Mumford type in general. In other words, Andr\'e's $p$-adic
arithmetic triangle groups do not always come out with the
uniformization by the spaces of form 
$\P^{1,\an}_{\C_p}-\limitset_{\Gamma}$, but by \'etale coverings 
of them. 
The following arithmetic triangle groups do not appear in our list:

\begin{condlist}
\item [\ \ $p=5$:] $(2,5,10)$, $(5,5,5)$, $(2,15,30)$, $(3,10,30)$,
      $(15,15,15)$. 
\item [\ \ $p=3$:] $(3,6,6)$, $(6,12,12)$, $(9,18,18)$.
\item [\ \ $p=2$:] $(2,4,8)$, $(2,8,8)$, $(4,4,4)$, $(4,8,8)$, $(3,4,12)$,
$(2,8,16)$, $(4,16,16)$, $(8,8,8)$, $(2,12,24)$, $(3,8,24)$,
$(6,24,24)$, $(12,12,12)$.
\end{condlist}

\paragraph\label{para-outline}\
{\sl Outline of the proof.}\ 
The proof of the theorem will be carried out by studying the action on
(a subtree of) the Bruhat-Tits tree (cf.\
\cite{Mum72}\cite{GvP80}\cite{Ser80}) by discrete subgroups in
$\PGL(2,\C_p)$.
To a finitely generated discrete subgroup $\Gamma$, there exists
associated tree $\tree_{\Gamma}$ (resp.\ $\tree^{\ast}_{\Gamma}$), which is
the tree generated by the set of limit points (resp.\ limit points
together with fixed points of elliptic elements) in
$\P^{1,\mathrm{an}}_{\C_p}$, being as the set of ends
(see \ref{para-review}).
The tree $\tree_{\Gamma}$ is equal to the Mumford's tree
(e.g.\ if $\Gamma$ is a Schottky subgroup), but the other tree 
is in general ``bigger''.
The advantage of the tree $\tree^{\ast}_{\Gamma}$
is its link with the ramification (or, branch) points; more precisely,
there exists a canonical bijection between branch points of the
uniformization map
$\varpi_{\Gamma}\colon(\P^{1,\mathrm{an}}_{\C_p}-\limitset_{\Gamma})
\rightarrow\Gamma\backslash(\P^{1,\mathrm{an}}_{\C_p}-\limitset_{\Gamma})$
and ends of the quotient graph
$T^{\ast}_{\Gamma}=\Gamma\backslash\tree^{\ast}_{\Gamma}$ (Proposition 
\ref{pro-branch}).
By this, we have the following principle:

\vspace{1ex}
{\bf Proposition.}\ 
{\sl A finitely generated discrete subgroup $\Gamma\subset\PGL(2,\C_p)$
is a $p$-adic triangle group of Mumford type if and only if the graph
$T^{\ast}_{\Gamma}$ is a tree having precisely three ends.
}\qed

\vspace{1ex}
That the graph $T^{\ast}_{\Gamma}$ is a tree is equivalent to 
that the quotient
$\Gamma\backslash(\P^{1,\mathrm{an}}_{\C_p}-\limitset_{\Gamma})$ is a
curve of genus $0$.
The formation of the trees $\tree^{\ast}_{\Gamma}$ admits the following obvious
functoriality: For an inclusion $\Gamma_1\subseteq\Gamma_2$ of finitely
generated subgroups, we have an inclusion of trees $\tree^{\ast}_{\Gamma_1}
\subseteq\tree^{\ast}_{\Gamma_2}$, and hence, a map
$T^{\ast}_{\Gamma_1}\rightarrow T^{\ast}_{\Gamma_2}$.

We decorate the tree $T^{\ast}_{\Gamma}$ by groups attached to
vertices and edges, which are simply the stabilizers of them. This gives 
rise to the data, so called, the {\it tree of groups}
$(T^{\ast}_{\Gamma},\Gamma_{\bullet})$.
The main point of the proof is that, essentially, the data
$(T^{\ast}_{\Gamma},\Gamma_{\bullet})$, considered as abstract tree
of groups, recovers $\Gamma$. Needless to say, to recover
$\Gamma$, one 
has to take nice embeddings of groups $\Gamma_{\bullet}$ in $\PGL(2,K)$. 
This has been discuss in \cite{Kat00}, where a complete criterion
for an abstract tree of groups $(T,G_{\bullet})$ to be realizable was
given. 
By this, the task in the proof reduces basically to
a purely combinatorial problem: Classify all possible $\ast$-admissible
tree of groups with exactly three ends.
This combinatorial business is easy in principle, but requires a lot of
care. We will introduce a notion of push-out, or direct limit of trees
of groups, which will be helpful to carry out the combinatorics.
In proving the theorem, we will exhibit
all the trees of groups in simple pictures, by which, besides, the
abstract group structure of the corresponding $\Gamma$ can be deduced. 

\vspace{1ex}
\paragraph\label{para-notation}\
{\sl Notation and conventions}.\ 
Throughout this paper $K$ denotes a finite extension of $\Q_p$,
$\O_K$ the integer ring, and $\pi$ a prime element in $\O_K$. 
We write $[K\colon\Q_p]=ef$, where $e$ is the ramification degree and
$q=p^f$ is the the number of elements in the residue field
$k=\O_K/\pi\O_K$.
We denote by $\valuation\colon K^{\times}\rightarrow\Z$ 
the normalized (i.e., $\valuation(\pi)=1$) valuation. 

For an abstract tree $T$ we denote by $\Vert(T)$ (resp.\
$\Edge(T)$, $\Ends(T)$) the set of all vertices (resp.\
unoriented edges, ends).
The notation $v\vdash\sigma$ for $v\in\Vert(T)$ and
$\sigma\in\Edge(T)$ means that $\sigma$ emanates from $v$. 
For a vertex $v\in\Vert(T)$ we denote by $\Star_v(T)$ the set of 
edges in $\Edge(T)$ emenating from $v$.
For two vertices $v_0$ and $v_1$, we denote by $[v_0,v_1]$ the geodesic
path connecting them.
For $\varepsilon_0,\varepsilon_1\in\Ends(T)$ and $v\in\Vert(T)$,
the unique straight-line (resp.\ half-line) connecting $\varepsilon_0$
and $\varepsilon_1$ (resp.\ $v$ and $\varepsilon_0$) is denoted by
$]\varepsilon_0,\varepsilon_1[$ (resp.\ $[v,\varepsilon_0[$).
The geometric realization $|T|$ is metrized so that the path
$[v_0,v_1]$ ($v_0,v_1\in\Vert(T)$) is of length equal to the number
of edges in it. The metric function is denoted by
$d_T(\cdot,\cdot)$, or simply by $d(\cdot,\cdot)$.
If $T$ is a subtree of $\tree_K$, the Bruhat-Tits tree attached to
$\PGL(2,K)$, then we always regard the set
$\Ends(T)$ as a subset of $\P^1(K)$.
In dealing with a tree we often switch to regard it as a topological
space by means of the geometric realization.

\vspace{2ex}
\sectioning\label{section-preliminaries}
{\bf Preliminaries}

\vspace{1ex}
Let us first review the basic facts on trees and groups which were dealt
with in \cite{Kat00}:

\vspace{1ex}
\paragraph\label{para-review}\
{\sl Trees and groups}\ (cf.\ \cite[\S 2]{Kat00}).
Let $\Gamma$ be a finitely generated discrete subgroup in $\PGL(2,K)$,
and suppose that $K$ has been chosen to be large enough such that the
fixed points (in $\P^1_K$) of any elliptic element in $\Gamma$ are
$K$-valued. Such a $\Gamma$ associates a subtree $\tree^{\ast}_{\Gamma}$ 
in the Bruhat-Tits tree $\tree_K$ characterized by (i) the set of ends
of $\tree^{\ast}_{\Gamma}$ are in the canonical bijection with the
closure of the set of fixed points of each 
element ($\neq 1$) in $\Gamma$, and (ii) $\tree^{\ast}_{\Gamma}$ is
minimal among subtrees having this property. Clearly,
$\tree^{\ast}_{\Gamma}$ is acted on by $\Gamma$. 
Attaching the stabilizers to each vertex and edge, we get the tree of
groups $(\tree^{\ast}_{\Gamma},\Gamma_{\bullet})$.
Suppose that the quotient
$T^{\ast}_{\Gamma}=\Gamma\backslash\tree^{\ast}_{\Gamma}$ is a tree,
by which one can consider a section $\iota_{\Gamma}\colon
T^{\ast}_{\Gamma}\hookrightarrow\tree^{\ast}_{\Gamma}$.
The section $\iota_{\Gamma}$ gives rise to a tree of groups
$(T^{\ast}_{\Gamma},\Gamma_{\bullet})$ in the obvious way.

\vspace{1ex}
{\sl Convention.}\ Throughout this paper, when a finitely generated
discrete subgroup $\Gamma$ is discussed, the field $K$ is assumed to be
chosen large enough so that the tree $\tree^{\ast}_K$ can be defined.

\vspace{1ex}
\paragraph\label{pro-branch}\
{\bf Proposition.}\ 
{\sl There exist canonical bijections, compatible with the quotient
maps, 
$$
\begin{array}{ccc}
\left\{
\begin{minipage}{4cm}
\setlength{\baselineskip}{.85\baselineskip}
\begin{small}
{\slshape Ramification points in $\Omega_{\Gamma}$ of the map
$\varpi_{\Gamma}$}
\end{small}
\end{minipage}
\right\}
&\longleftrightarrow&
\Ends(\tree^{\ast}_{\Gamma})-\Ends(\tree_{\Gamma})
\\
\bigdownarrow&&\bigdownarrow\\
\left\{
\begin{minipage}{4cm}
\setlength{\baselineskip}{.85\baselineskip}
\begin{small}
{\slshape Branch points in $\Gamma\backslash\Omega_{\Gamma}$}
\end{small}
\end{minipage}
\right\}
&\longleftrightarrow&
\Ends(T^{\ast}_{\Gamma})
\rlap{.}
\end{array}
$$
Moreover, the decomposition group of a ramification point coincides with 
the stabilizer of the corresponding end.}
\qed

\vspace{1ex}
\paragraph\label{para-admissible}\
{\sl $\ast$-admissibility}\ (\cite[\S 3]{Kat00}).
Conversely, let $(T,G_{\bullet})$ be an abstract tree of groups, and
suppose that we 
are given embeddings $G_v\hookrightarrow\PGL(2,K)$ for any $v\in\Vert(T)$
compatible with each $G_{\sigma}\hookrightarrow G_v$ for any
$v\vdash\sigma$. Such embeddings 
gives rise to subtrees $\tree^{\ast}_{G_v}$ (for $v\in\Vert(T)$) as in
\ref{para-review}.
Let $\til{\tree}_{G_{\bullet}}$ be the minimal subtree in $\tree_K$
containing all $\tree^{\ast}_{G_v}$.
The set of ends in $\til{\tree}_{G_{\bullet}}$ is the union
of the set of ends in $\tree^{\ast}_{G_v}$ for $v\in\Vert(T)$.
This tree is labelled by groups (not necessarily finite)
$\til{G}_{\bullet}$ as follows: For a vertex
$v\in\Vert(\til{\tree}_{G_{\bullet}})$ the group $\til{G}_v$ is the
subgroup in $\PGL(2,K)$ generated by $(G_u)_v$ (the stabilizer at $v$ by 
the action of $G_u$ on $\tree_K$) for all $u\in\Vert(T)$; the definition 
of the group $\til{G}_{\sigma}$ for $\sigma\in\Edge(\til{\tree}_{G_{\bullet}})$
is similar, which is just the intersection of $\til{G}_v$'s at the two
extremities.

\vspace{1ex}
\paragraph\label{def-admissible}\
{\bf Definition.}\ An {\it admissible embedding} of an abstract tree of
groups $(T,G_{\bullet})$ is an embedding $\iota\colon T\hookrightarrow\tree_K$ 
of trees together with embeddings $G_v\hookrightarrow\PGL(2,K)$ for any
$v\in\Vert(T)$ compatible with each $G_{\sigma}\hookrightarrow G_v$ for
any $v\vdash\sigma$ such that the following conditions are satisfied:
\begin{condlist}
\item[(1)] $\iota(T)\subset\til{\tree}_{G_{\bullet}}$.
\item[(2)] For any $v\in\Vert(T)$ and $\gamma\in G_v$ ($\gamma\neq 1$),
           there exists $\delta\in\Gamma$ such that
           $M(\delta\gamma\delta^{-1})\bigcap\iota(T)$ contains an edge, 
           where $\Gamma$ is the subgroup in $\PGL(2,K)$ generated by
           all $G_v$ for $v\in\Vert(T)$.
\item[(3)] $\til{G}_{\iota(v)}=G_v$ for any $v\in\Vert(T)$. 
\item[(4)] $\til{G}_{\iota(\sigma)}=G_{\sigma}$ for any
           $\sigma\in\Edge(T)$.
\item[(5)] For any $v\in\Vert(T)$, we have $\Star_v(T)\cong G_v\backslash
           (G_v\cdot\Star_{\iota(v)}(\til{\tree}_{G_{\bullet}}))$ by the
           composite of $\iota$ followed by the projection.
\end{condlist}
An abstract tree of groups $(T,G_{\bullet})$ is said
to be {\it $\ast$-admissible} if it has an admissible embedding and the
associated amalgam $\lim_{\rightarrow}(T,G_{\bullet})$ is finitely
generated.

\vspace{1ex}
The tree of groups $(T^{\ast}_{\Gamma},\Gamma_{\bullet})$ associated to
a finitely generated discrete subgroup $\Gamma$ is $\ast$-admissible by
the section $\iota_{\Gamma}$ (cf.\ \cite[3.5]{Kat00}).
The following theorem, proved in \cite[\S 3]{Kat00}, states the
converse:

\vspace{1ex}
\paragraph\label{thm-realization}\
{\bf Theorem.}\ 
{\sl Let $(T,G_{\bullet})$ be a $\ast$-admissible tree of groups and 
$\iota\colon T\hookrightarrow\til{\tree}_{G_{\bullet}}$ with 
$\{G_v\hookrightarrow\PGL(2,K)\}_{v\in\Vert(T)}$ an admissible
embedding. Let $\Gamma$ be the subgroup in $\PGL(2,K)$ generated by all
$G_v$ for $v\in\Vert(T)$ and set
$$
\tree^{\ast}=\bigcup_{\gamma\in\Gamma}\gamma\cdot\iota(T)
$$
in $\tree_K$. Then:

\vspace{1ex}
(1) The group $\Gamma$ is a finitely generated discrete subgroup in
$\PGL(2,K)$ isomorphic to $\lim_{\rightarrow}(T,G_{\bullet})$.

\vspace{1ex}
(2) The subset $\tree^{\ast}$ in $\tree_K$ is a tree and
$\Gamma\backslash\tree^{\ast}\cong T$.

\vspace{1ex}
(3) The embedding $\iota$ gives a section $T\hookrightarrow\tree^{\ast}$ 
by which the induced tree of groups $(T,\Gamma_{\bullet})$ equals to 
$(T,G_{\bullet})$.

\vspace{1ex}
Moreover, if $\tree^{\ast}_{\Gamma}$ is the tree associated to $\Gamma$ as
in \ref{para-review}, then
$\tree^{\ast}_{\Gamma}\subseteq\tree^{\ast}$, and the induced inclusion
$T^{\ast}_{\Gamma}\hookrightarrow T$ enjoys the following:

\vspace{1ex}
(4) The induced inclusion
$\Ends(T^{\ast}_{\Gamma})\hookrightarrow\Ends(T)$ is a bijection.

\vspace{1ex}
(5) The section $\iota$ restricts to a section
$T^{\ast}_{\Gamma}\hookrightarrow\tree^{\ast}_{\Gamma}$ by which the
induced tree of groups $(T^{\ast}_{\Gamma},\Gamma_{\bullet})$ is the
restriction of $(T,\Gamma_{\bullet})=(T,G_{\bullet})$.

\vspace{1ex}
(6) The tree of groups $(T^{\ast}_{\Gamma},\Gamma_{\bullet})$ is a
contraction of $(T,\Gamma_{\bullet})=(T,G_{\bullet})$.}
\qed

\vspace{1ex}
\paragraph\label{lem-minimal}\
{\bf Lemma.}\ {\sl Let $\Gamma$ be a finitely generated discrete
subgroup such that
$T^{\ast}_{\Gamma}=\Gamma\backslash\tree^{\ast}_{\Gamma}$ is contraction 
minimal, i.e., there is no proper subtree in $T^{\ast}_{\Gamma}$ having
the same set of ends. }

\vspace{1ex}
\pf
Otherwise, there exists a vertex $v\in\Vert(T^{\ast}_{\Gamma})$ such
that $T-\{v\}$ is connected. Hence such a vertex occurs also in
$\tree^{\ast}_{\Gamma}$; but this contradicts that
$\tree^{\ast}_{\Gamma}$ is the smallest one having the prescribed set of 
ends.
\qed

\vspace{1ex}
\paragraph\label{def-tripod}\
{\bf Definition.}\ A {\it tripod} is a tree which is the union of three
half-lines $\ell_i$ $(i=1,2,3)$ starting at a common vertex $v$, called
the {\it center}, such that $\ell_i\bigcap\ell_j=\{v\}$ for any
$i\neq j$.

\vspace{1ex}
\paragraph\label{cor-minimal}\
{\bf Corollary.}\ {\sl If $\Gamma$ is a $p$-adic Schwarzian triangle
groups of Mumford type, then $T^{\ast}_{\Gamma}$ is a tripod.}

\vspace{1ex}
\pf
Clear from Proposition \ref{para-outline} and Lemma \ref{lem-minimal}.
\qed

\vspace{1ex}
\paragraph\label{para-mirror}
Finally, we recall the following: Let $\gamma$ $(\neq 1)$ be an
elliptic element of finite order in $\PGL(2,K)$, and suppose $K$ is
taken to be large 
enough for the fixed points of $\gamma$ to be $K$-valued. The apartment
connecting two fixed points of $\gamma$ is called the {\it mirror} of
$\gamma$, and denoted by $M(\gamma)$.
Needless to say, it is contained in the fixed locus by $\gamma$ in
$\tree_K$.

\vspace{1ex}
\paragraph\label{lem-fixedlocus}\
{\bf Lemma}\ (\cite[2.10]{Kat00}).\ 
{\sl Let $n$ be the order of $\gamma$, and set
$G=\langle\gamma\rangle$. 

(1) Let $v_0\in M(\gamma)$. If $(n,p)=1$, then $G$ acts
freely on the $q-1$ vertices adjacent to $v_0$ not lying on
$M(\gamma)$, where $q$ is the number of elements in the residue field
$k$. 

(2) Suppose $n=p^r$ for $r\geq 1$, and set
$s=\valuation(\zeta_{p^r}-1)$, where $\zeta_{p^r}$ is a primitive
$p^r$-th root of unity, and $\valuation$ is the normalized (i.e.\
$\valuation(\pi)=1$) valuation.
Then a vertex $v\in\tree_K$ is fixed by $G$
if and only if $0\leq d(v,M(\gamma))\leq s$.}
\qed

\vspace{2ex}
\sectioning\label{section-directlimit}
{\bf Direct limit of trees}

\vspace{1ex}
\paragraph\label{para-pushout}
Given a diagram
$\tree_1\leftarrow\tree_0\rightarrow\tree_2$ of morphisms of trees, one
can define the push-out 
$$
\begin{array}{ccc}
\tree_0&\longrightarrow&\tree_1\\
\bigdownarrow&&\bigdownarrow\\
\tree_2&\longrightarrow&\tree_1\rlap{$\sharp_{\stree_0}\tree_2$}
\end{array}
$$
in the category of graphs; it is, regarded as a diagram of topological
spaces, simply the push-out in the category of topological
spaces. In slightly more formal terms, the graph
$\tree_1\sharp_{\stree_0}\tree_2$ has the set of vertices
$\Vert(\tree_1)\coprod_{\Vert(\stree_0)}\Vert(\tree_2)$ (push-out of
sets) and the similarly defined set of oriented edges together with the
naturally defined notion of origin and terminus of edges. It is clear
that the push-out $\tree_1\sharp_{\stree_0}\tree_2$ of trees is again a
tree, provided that $\tree_0$ is not empty.

The similar construction can be applied for push-out of trees of groups, 
where a morphism
$\phi\colon(\tree_0,G_{0,\bullet})\rightarrow(\tree_1,G_{1,\bullet})$ of
trees of groups is defined to be a morphism of trees
$\phi\colon\tree_0\rightarrow\tree_1$ together with the collection of
monomorphisms of groups $G_{\ast}\rightarrow G_{\phi(\ast)}$.
For a diagram $(\tree_1,G_{1,\bullet})\leftarrow(\tree_0,G_{0,\bullet})
\rightarrow(\tree_2,G_{2,\bullet})$ of trees of groups, the push-out 
$(\tree_1\sharp_{\stree_0}\tree_2,G_{1_02,\bullet})$ is endowed with the 
amalgam groups; to see that it is actually a tree of groups, one has to
show that for $v\vdash\sigma$ in $\tree_1\sharp_{\stree_0}\tree_2$ the
induced morphism $G_{1_02,\sigma}\rightarrow G_{1_02,v}$ is
injective. This follows from the structure theorem of amalgam groups 
\cite[I.1.2]{Ser80}.

\vspace{1ex}
\paragraph\label{exa-pushout}\
{\bf Example.}\ Let $\tree_0$ be the straight-line 
$$
\tree_0=\cdots\textrm{---}v_{-2}\textrm{---}v_{-1}\textrm{---}v_0\textrm{---}
v_1\textrm{---}v_2\textrm{---}\cdots,
$$
and $\tree_1$ and $\tree_2$ half-lines
\begin{eqnarray*}
\tree_1&=&u_0\textrm{---}u_1\textrm{---}u_2\textrm{---}\cdots,\\
\tree_2&=&\cdots\textrm{---}w_{-2}\textrm{---}w_{-1}\textrm{---}w_0.\\
\end{eqnarray*}

\vspace{-2ex}\noindent
Let $m$ be a positive integer.
The morphisms $f\colon\tree_0\rightarrow\tree_1$ and
$g\colon\tree_0\rightarrow\tree_2$ are defined by $f(v_n)=u_{|n|}$ and 
$g(v_n)=w_{-|n-m|}$, respectively.
Then the push-out $\tree_1\sharp_{\stree_0}\tree_2$ is a segment of
length $m$, isomorphic to $[u_0,u_m]$ (and $[w_{-m},w_0]$).

\vspace{1ex}
\paragraph\label{para-construction}\ 
Let $\Gamma$ a finitely generated discrete subgroup (e.g.\ a finite
subgroup) in $\PGL(2,K)$, and
$G\subseteq\Gamma$ a finite subgroup. 
By the construction of $\ast$-trees (\ref{para-review}, cf.\ \cite[\S
2]{Kat00}), we have an inclusion of subtrees
$\tree^{\ast}_G\subseteq\tree^{\ast}_{\Gamma}$, which yields the commutative
diagram  
$$
\begin{array}{ccc}
\tree^{\ast}_G&\longhookrightarrow&\tree^{\ast}_{\Gamma}\\
\llap{$\scriptstyle{\varrho_G}$}\bigdownarrow&&\bigdownarrow
\rlap{$\scriptstyle{\varrho_{\Gamma}}$}\\
T^{\ast}_G&\underrel{\longrightarrow}{\varrho^G_{\Gamma}}&
T^{\ast}_{\Gamma}
\end{array}
$$
of graphs, where $\varrho_G$ and $\varrho_{\Gamma}$ are quotients by $G$
and $\Gamma$, respectively. 
The morphism $\varrho^G_{\Gamma}$ is not in general injective. 
Note that (as one can see in Appendix below) the quotient graph
$T^{\ast}_G$ for a finits subgroup is always a tree.
Suppose that $T^{\ast}_{\Gamma}$ is a tree.
The four trees in the above diagram are then endowed with finite groups
as in \ref{para-review}, and become trees of groups. The above diagram
has an obvious extension to a diagram of trees of groups, where the
morphism between attached groups are defined abstractly.

Now suppose that we are given two finite subgroup $G_1$ and $G_2$ in
$\PGL(2,K)$ with $G_0=G_1\bigcap G_2$ such that $G_0\neq\{1\}$
Then one can consider the push-out diagram
$$
\begin{array}{ccc}
T^{\ast}_{G_0}&\stackrel{\varrho^{G_0}_{G_1}}{\longrightarrow}&T^{\ast}_{G_1}\\
\llap{$\scriptstyle{\varrho^{G_0}_{G_2}}$}\bigdownarrow&&\bigdownarrow\\
T^{\ast}_{G_2}&\longrightarrow&
T^{\ast}_{G_1}\rlap{$\sharp_{T^{\ast}_{G_0}}T^{\ast}_{G_2}$}
\end{array}
$$
of abstract trees of groups.

\vspace{1ex}
\paragraph\label{lem-finite}\
{\bf Lemma.}\ {\sl The following conditions are equivalent:
\begin{condlist}
\item[(1)] For any 
           $v\in\Vert(T^{\ast}_{G_1}\sharp_{T^{\ast}_{G_0}}T^{\ast}_{G_2})$,
           the group attached to $v$ is finite.
\item[(2)] For any $v\in\Vert(\tree^{\ast}_{G_0})$, either 
           $G_{1,v}\subseteq G_{2,v}$ or $G_{2,v}\subseteq G_{1,v}$ holds.
\end{condlist}}

\pf
Clear by the construction of push-out of trees of groups and the
structure theorem of amalgam groups \cite[I.1.2]{Ser80}.
\qed

\vspace{1ex}
\paragraph\label{para-directlimit}\ 
Let $\Gamma$ be a finitely generated discrete subgroup in $\PGL(2,K)$
such that $T^{\ast}_{\Gamma}$ is a tree, and $\iota_{\Gamma}\colon
T^{\ast}_{\Gamma}\hookrightarrow\tree^{\ast}_{\Gamma}$ a section.
The construction as in \ref{para-construction} can be carried out for 
$\Gamma_{v_1}$ and $\Gamma_{v_2}$ with the intersection
$\Gamma_{\sigma}$ coming from 
each edge $\sigma$ with extremities $v_1$ and $v_2$ in
$T^{\ast}_{\Gamma}$ as far as $\Gamma_{\sigma}\neq\{1\}$.
Hence we can define the {\it direct limit} along $T^{\ast}_{\Gamma}$
(similarly as in \cite[I.1.2, below Theorem 2]{Ser80}), denoted by 
$$
\lim_{\longrightarrow}T^{\ast}_{\Gamma_{\bullet}},
$$
which is a disjoint union of trees of groups; it is a single tree of
group, if there is no edge in $T^{\ast}_{\Gamma}$ with trivial group. 
The morphisms $\varrho^{\Gamma_v}_{\Gamma}$
($v\in\Vert(T^{\ast}_{\Gamma})$) and $\varrho^{\Gamma_{\sigma}}_{\Gamma}$
($\sigma\in\Edge(T^{\ast}_{\Gamma})$) induce a morphism 
$$
\lim_{\longrightarrow}T^{\ast}_{\Gamma_{\bullet}}
\longrightarrow
T^{\ast}_{\Gamma}.
\leqno{\indent\textrm{(\ref{para-directlimit}.1)}}
$$
It is clear that 
the image $S$ of this morphism is the union of the images of
$\tree^{\ast}_{\Gamma_v}$ for all $v\in\Vert(T^{\ast}_{\Gamma})$ under the
quotient map $\varrho_{\Gamma}$.
Let $S=\coprod_{i\in I}S_i$ be the decomposition into connected
components, and $s_{ij}$ ($i,j\in I$, $i\neq j$) the geodesic
path connecting $S_i$ and $S_j$.

\vspace{1ex}
\paragraph\label{lem-reducible1}\
{\bf Lemma.}\ {\sl If $s_{ij}$ does not meet any $S_k$ for
$k\not\in\{i,j\}$, then it contains an edge $\sigma$ with
$\Gamma_{\sigma}=\{1\}$.}

\vspace{1ex}
\pf
By Lemma \ref{lem-fixedlocus}, there exist two increasing sequences of
$p$-subgroups consisting of subgroups in the stabilizers of vertices in
$s_{ij}$, which are increasing ordered 
approaching to each $S_i$ and $S_j$. If there is no $\sigma$ with 
$\Gamma_{\sigma}=\{1\}$, then there exists a vertex $v$ in $s_{ij}$
fixed by two non-trivial $p$-groups having distinct mirrors. 
If $v\neq v_1$ and $v\neq v_2$, then, $\Gamma_v$ is not
contained in $\Gamma_w$'s for any vertex $w$ in $S$, and hence, $v$ is in the
image of $\tree^{\ast}_{\Gamma_v}$, thereby the contradiction.
If $v\in S_i$, then $\tree^{\ast}_{\Gamma_v}$ contains a mirror which is
mapped to $S_j$, and hence $s_{ij}$ is in the image of
$\tree^{\ast}_{\Gamma_v}$. 
\qed

\vspace{1ex}
\paragraph\label{lem-reducible2}\
{\bf Lemma.}\ {\sl Let $(T,G_{\bullet})$ be a $\ast$-admissible tree of
groups, and $\sigma\in\Edge(T)$ an edge such that $G_{\sigma}=\{1\}$.
Decompose $T=T_1\bigcup[v_1,v_2]\bigcup T_2$, where $v_1$ and $v_2$ are
the extremities of $\sigma$, such that $T_i\bigcap[v_1,v_2]=\{v_i\}$ for 
$i=1,2$. Then each $(T_i,G_{\bullet})$ is $\ast$-admissible, and
$\Gamma\cong\Gamma_1\ast\Gamma_2$, where
$\Gamma_i=\lim_{\rightarrow}(T_i,G_{\bullet})$.}

\vspace{1ex}
\pf
Clearly, we have $\Gamma\cong\Gamma_1\ast\Gamma_2$. Consider an
admissible embedding of $(T,G_{\bullet})$, which is restricted to each
$T_i$. To show that 
$(T_1,G_{\bullet})$ is $\ast$-admissible, only the condition
(\ref{def-admissible}.2) calls for a verification.
Let $v\in\Vert(T_1)$ and $\gamma\in G_v-\{1\}$. Take $\delta\in\Gamma$
such that $M(\delta\gamma\delta^{-1})\bigcap T\neq\emptyset$.
If $M(\delta\gamma\delta^{-1})\bigcap T_2\neq\emptyset$, then there
exists $w\in\Vert(T_2)$ such that $\chi=\delta\gamma\delta^{-1}$ belongs 
to $G_w$ yielding a non-trivial relation between elements in $\Gamma_1$
and $\Gamma_2$. Hence $M(\delta\gamma\delta^{-1})\bigcap
T_1\neq\emptyset$, thereby the lemma.
\qed

\vspace{1ex}
\paragraph\label{def-irreducible}\
{\bf Definition.}\ A tree of groups $(T,G_{\bullet})$
is said to be {\it irreducible} if $T$ does not contain an edge to which 
the trivial group is attached.

\vspace{1ex}
Due to Lemma \ref{lem-reducible1} and the minimality of
$T^{\ast}_{\Gamma}$ (Lemma \ref{lem-minimal}), 
if $T^{\ast}_{\Gamma}$ is irreducible, then the map
(\ref{para-directlimit}.1) of trees is surjective, i.e.\
$S=T^{\ast}_{\Gamma}$.

\vspace{1ex}
\paragraph\label{pro-directlimit}\
{\bf Proposition.}\ 
{\sl If $(T^{\ast}_{\Gamma},\Gamma_{\bullet})$ is irreducible, then 
the morphism (\ref{para-directlimit}.1) is an isomorphism of trees of
groups. In general, it is injective, and maps every connected 
component of $\lim_{\longrightarrow}T^{\ast}_{G_{\bullet}}$
isomorphically onto a subtree of groups in $T^{\ast}_{\Gamma}$.}

\vspace{1ex}
Before the proof, we need:

\vspace{1ex}
\paragraph\label{def-regular}\
{\bf Definition.}\ Let $G_1$ and $G_2$ be subgroups in
$\PGL(2,K)$. Then we say that $G_1$ and $G_2$ are {\it in regular
position} if, for any $K$-split torus $T$ in 
$\PGL(2,K)$, $G_1\bigcap T\neq\{1\}$ and $G_2\bigcap
T\neq\{1\}$ imply $G_1\bigcap G_2\bigcap T\neq\{1\}$.

\vspace{1ex}
\paragraph\label{lem-intersection}\
{\bf Lemma.}\ {\sl Let $v_1,v_2\in\Vert(T^{\ast}_{\Gamma})$.
Then $\Gamma_{v_1}$ and $\Gamma_{v_2}$ are in regular position.}

\vspace{1ex}
\pf
Let $T$ be a split torus, and $\gamma_1\in\Gamma_{v_1}\bigcap T-\{1\}$ and 
$\gamma_2\in\Gamma_{v_2}\bigcap T-\{1\}$.
Let $s=[v_1,v_2]$. Then $\Gamma_{v_1}\bigcap\Gamma_{v_2}$ is the set of
elements in $\Gamma$ which fix $s$ pointwise.
If the mirror $M=M(\gamma_1)=M(\gamma_2)$ contains $s$, then
$\gamma_1\in\Gamma_{v_1}\bigcap\Gamma_{v_2}$. If not, we have two cases: First,
if $M$ does not meet the interior of $s$, then, exchanging indices if
necessary, we may assume that $v_1$ is nearer to $M$ than $v_2$.
Then $s$ is fixed by $\gamma_2$, and hence, $\gamma_2\in
\Gamma_{\sigma}=\Gamma_{v_1}\bigcap\Gamma_{v_2}$.
If $M$ meets $s$ at a interior vertex $v$, then, by Lemma
\ref{lem-fixedlocus}, $\gamma_1$ and $\gamma_2$ are $p$-elements, and hence, 
exchanging indices if necessary, we may assume
$\langle\gamma_1\rangle\subseteq \langle\gamma_2\rangle$.
In this case, $s$ is fixed by $\gamma_2$.
\qed

\vspace{1ex}
For a finitely generated discrete subgroup $\Gamma$, $\fixset_{\Gamma}$
denotes the set of points in $\P^1_K$ fixed by an element in
$\Gamma-\{1\}$. Due to our convension about the field $K$, it consists 
of $K$-valued points.

\vspace{1ex}
\paragraph\label{lem-regular}\
{\bf Lemma.}\ {\sl Let $\Gamma_1,\Gamma_2\subset\PGL(2,K)$ be finitely
generated discrete subgroups which are in regular position.
Then we have
$\fixset_{\Gamma_1\bigcap\Gamma_2}=\fixset_{\Gamma_1}\bigcap
\fixset_{\Gamma_2}$.
In particular, if $\Gamma_1$ and $\Gamma_2$ are finite, then we have 
$\Ends(\tree^{\ast}_{\Gamma_1\bigcap\Gamma_2})=
\Ends(\tree^{\ast}_{\Gamma_1}\bigcap\tree^{\ast}_{\Gamma_2})$.}

\vspace{1ex}
\pf
$\fixset_{\Gamma_1\bigcap\Gamma_2}\subseteq\fixset_{\Gamma_1}\bigcap
\fixset_{\Gamma_2}$ is clear.
Let $z\in\fixset_{\Gamma_1}\bigcap\fixset_{\Gamma_2}$.
There exist $\gamma_1\in\Gamma_1-\{1\}$ and $\gamma_2\in\Gamma_2-\{1\}$
such that $\gamma_1(z)=\gamma_2(z)=z$.
Since no two element in a discrete subgroup share exactly one fixed
point (well-known, cf.\ \cite[2.5]{Kat00}), $\gamma_1$ and $\gamma_2$
belong to a same split torus. By the assumption, there exists $\gamma_3
\in\Gamma_1\bigcap\Gamma_2-\{1\}$ having the same fixed points as
$\gamma_1$ and $\gamma_2$, thereby the lemma.
\qed

\vspace{1ex}
{\it Proof of Proposition \ref{pro-directlimit}.}\ 
We embedd $T^{\ast}_{\Gamma}$ into $\tree^{\ast}_{\Gamma}$ (together
with attached groups) by the section $\iota_{\Gamma}$ fixed at the
beginning of our construction.
Let $S$ be the intersection of $T^{\ast}_{\Gamma}$ and the union of all
$\tree^{\ast}_{\Gamma_v}$ for $v\in\Vert(T^{\ast}_{\Gamma})$, which coincides 
with the image of the union of all $\tree^{\ast}_{\Gamma_v}$ under the
quotient map. We attach groups to $S$ in the obvious way. 
For any $v\in\Vert(T^{\ast}_{\Gamma})$, consider the natual morphism
$T^{\ast}_{\Gamma}\bigcap\tree^{\ast}_{\Gamma_v}\rightarrow
\lim_{\longrightarrow}T^{\ast}_{\Gamma_{\bullet}}$.
We can glue thus obtained morphisms to a morphism (together with
morphisms of groups) defined on
$S$; indeed, for $v_1,v_2\in\Vert(T^{\ast}_{\Gamma})$, 
Lemma \ref{lem-intersection}, Lemma \ref{lem-regular}, and Lemma
\ref{lem-minimal} imply that
$T^{\ast}_{\Gamma}\bigcap\tree^{\ast}_{\Gamma_{v_1}}\bigcap
\tree^{\ast}_{\Gamma_{v_2}}= 
T^{\ast}_{\Gamma}\bigcap\tree^{\ast}_{\Gamma_{v_1}\bigcap\Gamma_{v_2}}$.
It is easily verified that this morphism on $S$ gives the inverse of
(\ref{para-directlimit}.1). \qed

\vspace{2ex}
\sectioning\label{section-construction}
{\bf Construction of triangle groups}

\vspace{1ex}
In this section, we will construct $p$-adic Schwarzian triangle groups
of Mumford type in $p=2,3,5$. It will be proved in the next section that 
these are actually the only possible triangle groups.

\vspace{1ex}
\paragraph\label{para-5-1}\ {\sl $p=5$.}\ 
First we discuss in $p=5$, i.e., $K$ is a finite extention of $\Q_5$.
We begin with a finite subgroup $G_1\subset\PGL(2,K)$ isomorphic to
$A_5$. Inside $G_1$ we consider a subgroup $G_0\subset G_1$ isomorphic
to $D_5$. The morphism $\varrho^0_1\colon T^{\ast}_{G_0}\rightarrow
T^{\ast}_{G_1}$ is described as follows (see \ref{para-dihedral} and 
\ref{para-icosahedral} in Appendix; the pictures are
drawn obeying the convention in \ref{para-convention}):
$$
\begin{array}{c}
\setlength{\unitlength}{.7pt}
\begin{picture}(150,140)(0,0)
\put(81,52){\vector(2,-1){58}}
\put(81,52){\vector(2,-1){3}}
\put(75,55){\vector(-2,-1){64}}
\put(75,55){\circle*{4}}
\put(75,55){\line(0,1){18}}
\put(75,73){\vector(-1,3){18}}
\put(75,73){\circle*{4}}
\put(60,71){$\scriptscriptstyle{v_0}$}
\put(71,43){$\scriptscriptstyle{D_5}$}
\put(53,132){$\scriptscriptstyle{2}$}
\put(40,120){$\scriptscriptstyle{\varepsilon_0}$}
\put(3,17){$\scriptscriptstyle{2}$}
\put(6,32){$\scriptscriptstyle{\varepsilon_1}$}
\put(143,17){$\scriptscriptstyle{5}$}
\put(10,120){$\scriptstyle{T^{\ast}_{G_0}}$}
\end{picture}
\setlength{\unitlength}{.7pt}
\begin{picture}(100,140)(0,0)
\put(25,75){\vector(1,0){50}}
\put(45,83){$\scriptstyle{\varrho^0_1}$}
\end{picture}
\setlength{\unitlength}{.7pt}
\begin{picture}(150,140)(0,0)
\put(79,53){\vector(2,-1){60}}
\put(79,53){\vector(2,-1){3}}
\put(75,55){\vector(-2,-1){64}}
\put(75,55){\circle*{4}}
\put(75,73){\circle*{4}}
\put(75,80){\vector(0,1){47}}
\put(75,80){\vector(0,1){3}}
\put(75,66){\line(0,-1){11}}
\put(75,66){\vector(0,-1){3}}
\put(59,47){\circle*{4}}
\put(58,70){$\scriptscriptstyle{v'_0}$}
\put(80,57){$\scriptscriptstyle{v'_1}$}
\put(71,43){$\scriptscriptstyle{D_5}$}
\put(80,70){$\scriptscriptstyle{A_5}$}
\put(73,132){$\scriptscriptstyle{3}$}
\put(3,17){$\scriptscriptstyle{2}$}
\put(6,32){$\scriptscriptstyle{\varepsilon'_1}$}
\put(143,17){$\scriptscriptstyle{5}$}
\put(137,32){$\scriptscriptstyle{\varepsilon'_{\infty}}$}
\put(48,52){$\scriptscriptstyle{D_5}$}
\put(10,120){$\scriptstyle{T^{\ast}_{G_1}}$}
\end{picture}\\
\textrm{Figure 1}
\end{array}
$$
Here, $\varrho^0_1$ maps the locus below $v_0$ into
the locus below $v'_0$. The straight-line $]\varepsilon_0,\varepsilon_1[$
is mapped to the half-line $[v'_0,\varepsilon'_1[$ by a map like
$x\mapsto|x|$ with the folding at $v_0$.

Next we consider a finite subgroup $G_2$ isomorphic to $D_{10m}$ ($m\geq 
1$) such that $G_0\subseteq G_2$; $G_2$ is generated by $G_0$ and an
element of order $10m$ which commutes with elements of order $5$ in
$G_0$.
The morphism $\varrho^0_2\colon T^{\ast}_{G_0}\rightarrow
T^{\ast}_{G_2}$ is described as follows:
$$
\begin{array}{c}
\setlength{\unitlength}{.7pt}
\begin{picture}(150,140)(0,0)
\put(81,52){\vector(2,-1){58}}
\put(81,52){\vector(2,-1){3}}
\put(75,55){\vector(-2,-1){64}}
\put(75,55){\circle*{4}}
\put(75,55){\line(0,1){18}}
\put(75,73){\vector(-1,3){18}}
\put(75,73){\circle*{4}}
\put(60,59){$\scriptscriptstyle{v_1}$}
\put(71,43){$\scriptscriptstyle{D_5}$}
\put(53,132){$\scriptscriptstyle{2}$}
\put(40,120){$\scriptscriptstyle{\varepsilon_0}$}
\put(3,17){$\scriptscriptstyle{2}$}
\put(6,32){$\scriptscriptstyle{\varepsilon_1}$}
\put(143,17){$\scriptscriptstyle{5}$}
\put(10,120){$\scriptstyle{T^{\ast}_{G_0}}$}
\end{picture}
\setlength{\unitlength}{.7pt}
\begin{picture}(100,140)(0,0)
\put(25,75){\vector(1,0){50}}
\put(45,83){$\scriptstyle{\varrho^0_2}$}
\end{picture}
\setlength{\unitlength}{.7pt}
\begin{picture}(150,140)(0,0)
\put(81,52){\vector(2,-1){58}}
\put(81,52){\vector(2,-1){3}}
\put(69,52){\vector(-2,-1){58}}
\put(69,52){\vector(-2,-1){3}}
\put(75,55){\circle*{4}}
\put(75,62){\line(0,1){11}}
\put(75,62){\vector(0,1){3}}
\put(75,73){\vector(-1,3){18}}
\put(75,73){\circle*{4}}
\put(60,59){$\scriptscriptstyle{v''_1}$}
\put(67,40){$\scriptscriptstyle{D_{10m}}$}
\put(53,132){$\scriptscriptstyle{2}$}
\put(40,120){$\scriptscriptstyle{\varepsilon''_0}$}
\put(3,17){$\scriptscriptstyle{2}$}
\put(6,32){$\scriptscriptstyle{\varepsilon''_1}$}
\put(143,17){$\scriptscriptstyle{10m}$}
\put(10,120){$\scriptstyle{T^{\ast}_{G_2}}$}
\end{picture}\\
\textrm{Figure 2}
\end{array}
$$
Here the straight line $]\varepsilon_0,\varepsilon_1[$ is mapped to the
half-line $[v''_1,\varepsilon''_0[$.
By means of these data, it is not difficult to compute the push-out
$T=T^{\ast}_{G_1}\sharp_{T^{\ast}_{G_0}}T^{\ast}_{G_2}$; the only point
to pay attention is that the straight line
$]\varepsilon_0,\varepsilon_1[$ is mapped in $T$ onto a segment 
isomorphic to $[v_0,v_1]$ by the same reasoning as in Example
\ref{exa-pushout}. 
The resulting tree of groups $T$ looks like as in Figure 3.
$$
\begin{array}{c}
\setlength{\unitlength}{.7pt}
\begin{picture}(150,140)(0,0)
\put(81,52){\vector(2,-1){58}}
\put(81,52){\vector(2,-1){3}}
\put(69,52){\vector(-2,-1){58}}
\put(69,52){\vector(-2,-1){3}}
\put(75,55){\circle*{4}}
\put(75,73){\circle*{4}}
\put(75,80){\vector(0,1){47}}
\put(75,80){\vector(0,1){3}}
\put(75,67){\line(0,-1){6}}
\put(75,66){\vector(0,-1){3}}
\put(75,61){\vector(0,1){3}}
\put(59,47){\circle*{4}}
\put(67,40){$\scriptscriptstyle{D_{10m}}$}
\put(80,70){$\scriptscriptstyle{A_5}$}
\put(73,132){$\scriptscriptstyle{3}$}
\put(3,17){$\scriptscriptstyle{2}$}
\put(143,17){$\scriptscriptstyle{10m}$}
\put(48,52){$\scriptscriptstyle{D_5}$}
\put(-10,120){$\scriptstyle{
T^{\ast}_{G_1}\sharp_{T^{\ast}_{G_0}}T^{\ast}_{G_2}}$} 
\end{picture}\\
\textrm{Figure 3}
\end{array}
$$

\paragraph\label{para-5-2}
We are going to show that the tree of groups
$T=T^{\ast}_{G_1}\sharp_{T^{\ast}_{G_0}}T^{\ast}_{G_2}$ is
$\ast$-admissible. Let $\iota_{G_1}\colon
T^{\ast}_{G_1}\rightarrow\tree^{\ast}_{G_1}$ be a section.
Then we can find a section $\iota_{G_2}\colon
T^{\ast}_{G_2}\rightarrow\tree^{\ast}_{G_2}$ such that 
$\iota_{G_1}(T^{\ast}_{G_1})\bigcap\iota_{G_2}(T^{\ast}_{G_2})$
is the half-line
$[\iota_{G_1}(v'_0),\iota_{G_1}(\varepsilon'_{\infty})[$. 
Indeed, $\iota_{G_2}(T^{\ast}_{G_2})$ is the union of the following
three half-lines: 
(i) $[\iota_{G_1}(v'_1),\iota_{G_1}(\varepsilon'_{\infty})[$, 
(ii) a half-line in a mirror of an element of order $2$ in $G_1$
containing $[\iota_{G_1}(v'_0),\iota_{G_1}(v'_1),]$ and starting at
$\iota_{G_1}(v'_1)$, and (iii) a half-line starting at
$\iota_{G_1}(v'_1)$ contained in a mirror of an element of order $2$ in
$G_2$ not in $G_0$ and $Z_{10m}\subset G_2$.
Let us define $T^{\ast}_{G_1}\rightarrow\iota_{G_1}(T^{\ast}_{G_1})$
(resp.\ $T^{\ast}_{G_2}\rightarrow\iota_{G_2}(T^{\ast}_{G_2})$) which
coincides with $\iota_{G_1}$ (resp.\ $\iota_{G_2}$) except for that the
half-line $[v'_0,\varepsilon'_1]$ (resp.\ $[v''_1,\varepsilon''_0]$) is
mapped to $[\iota_{G_1}(v'_0),\iota_{G_1}(v'_1)]$ in an obvious way. 
This induces an embedding $T\rightarrow\tree_K$ of the tree $T$.
Then, together with $G_1,G_2\hookrightarrow\PGL(2,K)$, this gives an
admissible embedding; indeed, the local structure of $T$ at any vertex
is isomorphic to that around a vertex either in $T^{\ast}_{G_1}$ or in 
$T^{\ast}_{G_2}$, because our situation satisfies the condition (2) in 
Lemma \ref{lem-finite}. Hence (\ref{def-admissible}.3),
(\ref{def-admissible}.4), and (\ref{def-admissible}.5) are valid.
The validity of (\ref{def-admissible}.2) is evident, since $T$ contains
the images of all the mirrors in $T^{\ast}_{G_1}$ and $T^{\ast}_{G_2}$.

As a result, we get a triangle group of index $(2,3,10m)$ ($m\geq 1$),
isomorphic to the amalgam product $A_5\ast_{D_5}D_{10m}$.

\vspace{1ex}
\paragraph\label{para-5-3}
We can also find the following two $\ast$-admissible trees of groups:
$$
\begin{array}{c}
\setlength{\unitlength}{.7pt}
\begin{picture}(150,140)(0,0)
\put(79,53){\vector(2,-1){60}}
\put(79,53){\vector(2,-1){3}}
\put(75,55){\vector(-2,-1){64}}
\put(75,55){\circle*{4}}
\put(75,73){\circle*{4}}
\put(75,80){\vector(0,1){47}}
\put(75,80){\vector(0,1){3}}
\put(75,66){\line(0,-1){11}}
\put(75,66){\vector(0,-1){3}}
\put(59,47){\circle*{4}}
\put(79,57){$\scriptscriptstyle{D_{10m+5}}$}
\put(80,70){$\scriptscriptstyle{A_5}$}
\put(73,132){$\scriptscriptstyle{3}$}
\put(3,17){$\scriptscriptstyle{2}$}
\put(143,17){$\scriptscriptstyle{10m+5}$}
\put(48,52){$\scriptscriptstyle{D_5}$}
\end{picture}
\qquad\qquad
\setlength{\unitlength}{.7pt}
\begin{picture}(150,140)(0,0)
\put(75,55){\circle*{4}}
\put(75,62){\vector(0,1){65}}
\put(75,62){\vector(0,1){3}}
\put(75,55){\line(2,-1){12}}
\put(87,49){\vector(-2,1){3}}
\put(95,45){\vector(2,-1){44}}
\put(95,45){\vector(2,-1){3}}
\put(75,55){\line(-2,-1){12}}
\put(63,49){\vector(2,1){3}}
\put(55,45){\vector(-2,-1){44}}
\put(55,45){\vector(-2,-1){3}}
\put(91,47){\circle*{4}}
\put(59,47){\circle*{4}}
\put(69,43){$\scriptscriptstyle{D_{5l}}$}
\put(71,132){$\scriptscriptstyle{5l}$}
\put(3,17){$\scriptscriptstyle{3}$}
\put(143,17){$\scriptscriptstyle{3}$}
\put(92,52){$\scriptscriptstyle{A_5}$}
\put(48,52){$\scriptscriptstyle{A_5}$}
\end{picture}\\
\textrm{Figure 4}
\end{array}
$$
We sketch the construction of these trees, and details are left to the
reader: 
The construction of the first one in Figure 4 is similar to that of
Figure 3 as above; it is even simpler, because, in this case, the
morphism $\varrho^0_2\colon T^{\ast}_{D_5}\rightarrow
T^{\ast}_{D_{10m+5}}$ is injective.
We call in general the procedure like this a {\it replacement}; it is, so
to speak, the replacement of the mirror of order $5$ by that of order
$10m+5$. 
The second one with $l=1$ is by 
$T^{\ast}_{G_1}\sharp_{T^{\ast}_{G_0}}T^{\ast}_{G'_1}$, where $G'_1$ is
another embedded $A_5$ which is the twist of $G_1$ by the non-trivial
element in $N(D_5)/(N(D_5)\bigcap N(A_5))\cong Z_2$.
Here $N(G)$ for a subgroup $G\subseteq\PGL(2,K)$ stands for the
normalizer of $G$ in $\PGL(2,K)$; note that $N(D_5)=D_{10}$ and
$N(A_5)=A_5$.

The second one with $l$ odd is constructed by the method similar to that
of the first one (replacement of the mirror of order $5$ by 
that of order $10m+5$), started by the one with $l=1$.
The construction of the second one with $l$ even is outlined as follows: 
We start at the tree of groups $T$ as in Fugure 3. Let $\delta$ be the
involution in $D_{10m}$ not contained in $G_0$ nor in $Z_{10m}\subset
D_{10m}$.
Let $G'_1=\delta G_1\delta$ and $G'_0=\delta G_0\delta$.
(Note that the ``$D_5$'' denoted in Figure 3 is $G'_0$.)
Then we consider the push-out
$T\sharp_{T^{\ast}_{G'_0}}T^{\ast}_{G'_1}$, which can be described by a
similar method as in \ref{para-5-1}, and this gives the desired tree of
groups. 

These are shown to be $\ast$-admissible by means of appropriate
embeddings, constructed by an idea similar to that in \ref{para-5-2}.

\vspace{1ex}
In the following (until the end of this section), we will perform only 
sketchy constructions, but will present necessary data by which the
reader can verify the details at each step; all the trees presented
below are proved to be $\ast$-admissible by an appropriate embeddings
which can be easily found, similarly as above.

\vspace{1ex}
\paragraph\label{para-3-1}\ {\sl $p=3$.}\ 
The arguement similar to that in $p=5$, where $(A_5,D_5)$ is replaced by 
$(A_5,D_3)$ and $(S_4,D_3)$, works in $p=3$, which yields the following
six $\ast$-admissible trees of groups ($m,l\geq 1$)
$$
\begin{array}{c}
\setlength{\unitlength}{.7pt}
\begin{picture}(150,140)(0,0)
\put(81,52){\vector(2,-1){58}}
\put(81,52){\vector(2,-1){3}}
\put(69,52){\vector(-2,-1){58}}
\put(69,52){\vector(-2,-1){3}}
\put(75,55){\circle*{4}}
\put(75,73){\circle*{4}}
\put(75,80){\vector(0,1){47}}
\put(75,80){\vector(0,1){3}}
\put(75,67){\line(0,-1){6}}
\put(75,66){\vector(0,-1){3}}
\put(75,61){\vector(0,1){3}}
\put(59,47){\circle*{4}}
\put(67,40){$\scriptscriptstyle{D_{6m}}$}
\put(80,70){$\scriptscriptstyle{A_5}$}
\put(73,132){$\scriptscriptstyle{5}$}
\put(3,17){$\scriptscriptstyle{2}$}
\put(143,17){$\scriptscriptstyle{6m}$}
\put(48,52){$\scriptscriptstyle{D_3}$}
\end{picture}
\qquad
\setlength{\unitlength}{.7pt}
\begin{picture}(150,140)(0,0)
\put(79,53){\vector(2,-1){60}}
\put(79,53){\vector(2,-1){3}}
\put(75,55){\vector(-2,-1){64}}
\put(75,55){\circle*{4}}
\put(75,73){\circle*{4}}
\put(75,80){\vector(0,1){47}}
\put(75,80){\vector(0,1){3}}
\put(75,66){\line(0,-1){11}}
\put(75,66){\vector(0,-1){3}}
\put(59,47){\circle*{4}}
\put(79,57){$\scriptscriptstyle{D_{6m+3}}$}
\put(80,70){$\scriptscriptstyle{A_5}$}
\put(73,132){$\scriptscriptstyle{5}$}
\put(3,17){$\scriptscriptstyle{2}$}
\put(143,17){$\scriptscriptstyle{6m+3}$}
\put(48,52){$\scriptscriptstyle{D_3}$}
\end{picture}
\qquad
\setlength{\unitlength}{.7pt}
\begin{picture}(150,140)(0,0)
\put(75,55){\circle*{4}}
\put(75,62){\vector(0,1){65}}
\put(75,62){\vector(0,1){3}}
\put(75,55){\line(2,-1){12}}
\put(87,49){\vector(-2,1){3}}
\put(95,45){\vector(2,-1){44}}
\put(95,45){\vector(2,-1){3}}
\put(75,55){\line(-2,-1){12}}
\put(63,49){\vector(2,1){3}}
\put(55,45){\vector(-2,-1){44}}
\put(55,45){\vector(-2,-1){3}}
\put(91,47){\circle*{4}}
\put(59,47){\circle*{4}}
\put(69,43){$\scriptscriptstyle{D_{3l}}$}
\put(71,132){$\scriptscriptstyle{3l}$}
\put(3,17){$\scriptscriptstyle{5}$}
\put(143,17){$\scriptscriptstyle{5}$}
\put(92,52){$\scriptscriptstyle{A_5}$}
\put(48,52){$\scriptscriptstyle{A_5}$}
\end{picture}\\
\setlength{\unitlength}{.7pt}
\begin{picture}(150,140)(0,0)
\put(81,52){\vector(2,-1){58}}
\put(81,52){\vector(2,-1){3}}
\put(69,52){\vector(-2,-1){58}}
\put(69,52){\vector(-2,-1){3}}
\put(75,55){\circle*{4}}
\put(75,73){\circle*{4}}
\put(75,80){\vector(0,1){47}}
\put(75,80){\vector(0,1){3}}
\put(75,67){\line(0,-1){6}}
\put(75,66){\vector(0,-1){3}}
\put(75,61){\vector(0,1){3}}
\put(59,47){\circle*{4}}
\put(67,40){$\scriptscriptstyle{D_{6m}}$}
\put(80,70){$\scriptscriptstyle{S_4}$}
\put(73,132){$\scriptscriptstyle{4}$}
\put(3,17){$\scriptscriptstyle{2}$}
\put(143,17){$\scriptscriptstyle{6m}$}
\put(48,52){$\scriptscriptstyle{D_3}$}
\end{picture}
\qquad
\setlength{\unitlength}{.7pt}
\begin{picture}(150,140)(0,0)
\put(79,53){\vector(2,-1){60}}
\put(79,53){\vector(2,-1){3}}
\put(75,55){\vector(-2,-1){64}}
\put(75,55){\circle*{4}}
\put(75,73){\circle*{4}}
\put(75,80){\vector(0,1){47}}
\put(75,80){\vector(0,1){3}}
\put(75,66){\line(0,-1){11}}
\put(75,66){\vector(0,-1){3}}
\put(59,47){\circle*{4}}
\put(79,57){$\scriptscriptstyle{D_{6m+3}}$}
\put(80,70){$\scriptscriptstyle{S_4}$}
\put(73,132){$\scriptscriptstyle{4}$}
\put(3,17){$\scriptscriptstyle{2}$}
\put(143,17){$\scriptscriptstyle{6m+3}$}
\put(48,52){$\scriptscriptstyle{D_3}$}
\end{picture}
\qquad
\setlength{\unitlength}{.7pt}
\begin{picture}(150,140)(0,0)
\put(75,55){\circle*{4}}
\put(75,62){\vector(0,1){65}}
\put(75,62){\vector(0,1){3}}
\put(75,55){\line(2,-1){12}}
\put(87,49){\vector(-2,1){3}}
\put(95,45){\vector(2,-1){44}}
\put(95,45){\vector(2,-1){3}}
\put(75,55){\line(-2,-1){12}}
\put(63,49){\vector(2,1){3}}
\put(55,45){\vector(-2,-1){44}}
\put(55,45){\vector(-2,-1){3}}
\put(91,47){\circle*{4}}
\put(59,47){\circle*{4}}
\put(69,43){$\scriptscriptstyle{D_{3l}}$}
\put(71,132){$\scriptscriptstyle{3l}$}
\put(3,17){$\scriptscriptstyle{4}$}
\put(143,17){$\scriptscriptstyle{4}$}
\put(92,52){$\scriptscriptstyle{S_4}$}
\put(48,52){$\scriptscriptstyle{S_4}$}
\end{picture}\\
\textrm{Figure 5}
\end{array}
$$

Also, it is easy to find that the method to obtain the two trees in the
last column in Figure 5 is mixed up to get the one in Figure 6.
$$
\begin{array}{c}
\setlength{\unitlength}{.7pt}
\begin{picture}(150,140)(0,0)
\put(75,55){\circle*{4}}
\put(75,62){\vector(0,1){65}}
\put(75,62){\vector(0,1){3}}
\put(75,55){\line(2,-1){12}}
\put(87,49){\vector(-2,1){3}}
\put(95,45){\vector(2,-1){44}}
\put(95,45){\vector(2,-1){3}}
\put(75,55){\line(-2,-1){12}}
\put(63,49){\vector(2,1){3}}
\put(55,45){\vector(-2,-1){44}}
\put(55,45){\vector(-2,-1){3}}
\put(91,47){\circle*{4}}
\put(59,47){\circle*{4}}
\put(69,43){$\scriptscriptstyle{D_{3l}}$}
\put(71,132){$\scriptscriptstyle{3l}$}
\put(3,17){$\scriptscriptstyle{4}$}
\put(143,17){$\scriptscriptstyle{5}$}
\put(92,52){$\scriptscriptstyle{A_5}$}
\put(48,52){$\scriptscriptstyle{S_4}$}
\end{picture}\\
\textrm{Figure 6}
\end{array}
\qquad\qquad
\begin{array}{c}
\setlength{\unitlength}{.7pt}
\begin{picture}(150,140)(0,0)
\put(75,55){\vector(2,-1){64}}
\put(75,55){\vector(-2,-1){64}}
\put(75,55){\circle*{4}}
\put(75,73){\circle*{4}}
\multiput(75,56)(0,4){4}{\circle*{2}}
\put(75,80){\vector(0,1){47}}
\put(75,80){\vector(0,1){3}}
\put(72,43){$\scriptscriptstyle{3l}$}
\put(80,70){$\scriptscriptstyle{A_4}$}
\put(72,132){$\scriptscriptstyle{2}$}
\put(2,17){$\scriptscriptstyle{3l}$}
\put(142,17){$\scriptscriptstyle{3l}$}
\end{picture}\\
\textrm{Figure 7}
\end{array}
$$
The tree in Figure 7 can be found by considering
$T^{\ast}_{A_4}\sharp_{T^{\ast}_{Z_3}}T^{\ast}_{Z_{3l}}$, which is
nothing but the replacement of the mirror of order $3$ by that of order
$3l$.

\vspace{1ex}
\paragraph\label{para-2-1}\ {\sl $p=2$.}\ 
There are plenty of triangle groups in $p=2$. First, the push-out by
$(G_1,G_0,G_2)=(S_4,A_4,A_5)$ yields the first tree in Figure 8 with
$l=1$. 
$$
\begin{array}{c}
\setlength{\unitlength}{.7pt}
\begin{picture}(150,140)(0,0)
\put(75,55){\circle*{4}}
\put(75,73){\circle*{4}}
\put(75,109){\circle*{4}}
\multiput(75,55)(0,4){4}{\circle*{2}}
\put(75,78){\line(0,1){26}}
\put(75,78){\vector(0,1){3}}
\put(75,104){\vector(0,-1){3}}
\put(75,114){\vector(0,1){13}}
\put(75,114){\vector(0,1){3}}
\put(75,91){\circle*{4}}
\multiput(75,55)(4,-2){4}{\circle*{2}}
\put(95,45){\vector(2,-1){44}}
\put(95,45){\vector(2,-1){3}}
\put(91,47){\circle*{4}}
\multiput(75,55)(-4,-2){12}{\circle*{2}}
\put(59,47){\circle*{4}}
\put(43,39){\circle*{4}}
\put(27,31){\circle*{4}}
\put(23,29){\vector(-2,-1){12}}
\put(23,29){\vector(-2,-1){3}}
\put(59,47){\circle*{4}}
\put(80,70){$\scriptscriptstyle{S_4}$}
\put(80,88){$\scriptscriptstyle{A_4}$}
\put(80,106){$\scriptscriptstyle{A_5}$}
\put(73,132){$\scriptscriptstyle{5}$}
\put(2,17){$\scriptscriptstyle{2l}$}
\put(141,17){$\scriptscriptstyle{4m}$}
\put(91,52){$\scriptscriptstyle{D_{4m}}$}
\put(48,52){$\scriptscriptstyle{D_4}$}
\put(32,44){$\scriptscriptstyle{D_2}$}
\put(15,36){$\scriptscriptstyle{D_{2l}}$}
\end{picture}
\qquad
\setlength{\unitlength}{.7pt}
\begin{picture}(150,140)(0,0)
\put(75,55){\circle*{4}}
\put(75,73){\circle*{4}}
\multiput(75,56)(0,4){4}{\circle*{2}}
\put(75,80){\vector(0,1){47}}
\put(75,80){\vector(0,1){3}}
\multiput(87,49)(-4,2){4}{\circle*{2}}
\put(95,45){\line(2,-1){20}}
\put(95,45){\vector(2,-1){3}}
\put(115,35){\vector(-2,1){3}}
\put(123,31){\vector(2,-1){16}}
\put(123,31){\vector(2,-1){3}}
\put(105,40){\circle*{4}}
\put(119,33){\circle*{4}}
\multiput(63,49)(4,2){4}{\circle*{2}}
\put(55,45){\vector(-2,-1){44}}
\put(55,45){\vector(-2,-1){3}}
\put(45,40){\circle*{4}}
\put(31,33){\circle*{4}}
\put(91,47){\circle*{4}}
\put(59,47){\circle*{4}}
\put(78,70){$\scriptscriptstyle{D_{4m}}$}
\put(70,132){$\scriptscriptstyle{4m}$}
\put(3,17){$\scriptscriptstyle{3}$}
\put(143,17){$\scriptscriptstyle{5}$}
\put(92,52){$\scriptscriptstyle{S_4}$}
\put(106,45){$\scriptscriptstyle{A_4}$}
\put(120,38){$\scriptscriptstyle{A_5}$}
\put(48,52){$\scriptscriptstyle{S_4}$}
\put(34,45){$\scriptscriptstyle{A_4}$}
\put(20,38){$\scriptscriptstyle{A_4}$}
\end{picture}
\quad
\setlength{\unitlength}{.7pt}
\begin{picture}(150,140)(0,0)
\put(75,55){\circle*{4}}
\put(75,73){\circle*{4}}
\multiput(75,56)(0,4){4}{\circle*{2}}
\put(75,80){\vector(0,1){47}}
\put(75,80){\vector(0,1){3}}
\multiput(87,49)(-4,2){4}{\circle*{2}}
\put(95,45){\line(2,-1){20}}
\put(95,45){\vector(2,-1){3}}
\put(115,35){\vector(-2,1){3}}
\put(123,31){\vector(2,-1){16}}
\put(123,31){\vector(2,-1){3}}
\put(105,40){\circle*{4}}
\put(119,33){\circle*{4}}
\multiput(63,49)(4,2){4}{\circle*{2}}
\put(55,45){\line(-2,-1){20}}
\put(55,45){\vector(-2,-1){3}}
\put(35,35){\vector(2,1){3}}
\put(27,31){\vector(-2,-1){16}}
\put(27,31){\vector(-2,-1){3}}
\put(45,40){\circle*{4}}
\put(31,33){\circle*{4}}
\put(91,47){\circle*{4}}
\put(59,47){\circle*{4}}
\put(78,70){$\scriptscriptstyle{D_{4m}}$}
\put(70,132){$\scriptscriptstyle{4m}$}
\put(3,17){$\scriptscriptstyle{5}$}
\put(143,17){$\scriptscriptstyle{5}$}
\put(92,52){$\scriptscriptstyle{S_4}$}
\put(106,45){$\scriptscriptstyle{A_4}$}
\put(120,38){$\scriptscriptstyle{A_5}$}
\put(48,52){$\scriptscriptstyle{S_4}$}
\put(34,45){$\scriptscriptstyle{A_4}$}
\put(20,38){$\scriptscriptstyle{A_5}$}
\end{picture}
\\
\textrm{Figure 8}
\end{array}
$$
The third tree with $m=1$ is obtained from the first one $T$ (with $m=1$)
by the push-out $T\sharp_{T^{\ast}_{D_4}}T'$, 
where $T'$ is the twist $T'=\delta T$ by a non-trivial element in 
$N(D_4)/(N(D_4)\bigcap N(S_4))\cong Z_2$ (note that $N(D_4)=D_8$). 
The second one with $m=1$ is $T\sharp_{T^{\ast}_{D_4}}\delta
T^{\ast}_{G_1}$.
Replacing the mirror of order $4$ by that of order $4m$, we get the ones
in $m\geq 1$.
The center of the trees in Figure 8 is fixed by $D_4$ if
$m$ is odd, or by $D_8$ otherwise.
Note that the, in the last two trees in Figure 8, the mirror of order
$2$ of each $S_4$ is absorbed in the mirror of order $4m$, which is
easily seen by means of push-out. 

In the first tree in Figure 8, the mirror of order $2$ has been already
replaced by that of order $2l$; but we claim that this replacement
can be done if and only if $l$ is odd. 
To see this, let $G_1$ be an embedded $S_4$, and $G_0\subset G_1$ a
subgroup isomorphic to $D_2$.
The morphism $\varrho^0_1\colon T^{\ast}_{G_0}\rightarrow
T^{\ast}_{G_1}$ is described as follows:
$$
\begin{array}{c}
\setlength{\unitlength}{.7pt}
\begin{picture}(150,140)(0,0)
\put(75,55){\circle*{4}}
\put(75,73){\circle*{4}}
\multiput(75,56)(0,4){4}{\circle*{2}}
\put(75,80){\vector(0,1){47}}
\put(75,80){\vector(0,1){3}}
\multiput(75,55)(4,-2){4}{\circle*{2}}
\put(95,45){\vector(2,-1){44}}
\put(95,45){\vector(2,-1){3}}
\multiput(75,55)(-4,-2){4}{\circle*{2}}
\put(55,45){\vector(-2,-1){44}}
\put(55,45){\vector(-2,-1){3}}
\put(91,47){\circle*{4}}
\put(59,47){\circle*{4}}
\put(71,43){$\scriptscriptstyle{D_2}$}
\put(80,70){$\scriptscriptstyle{D_2}$}
\put(72,132){$\scriptscriptstyle{2}$}
\put(3,17){$\scriptscriptstyle{2}$}
\put(143,17){$\scriptscriptstyle{2}$}
\put(92,52){$\scriptscriptstyle{D_2}$}
\put(48,52){$\scriptscriptstyle{D_2}$}
\put(79,57){$\scriptscriptstyle{v_1}$}
\put(60,39){$\scriptscriptstyle{v_0}$}
\put(10,120){$\scriptstyle{T^{\ast}_{G_0}}$}
\end{picture}
\setlength{\unitlength}{.7pt}
\begin{picture}(100,140)(0,0)
\put(25,75){\vector(1,0){50}}
\put(45,83){$\scriptstyle{\varrho^0_1}$}
\end{picture}
\setlength{\unitlength}{.7pt}
\begin{picture}(150,140)(0,0)
\put(75,55){\circle*{4}}
\put(75,73){\circle*{4}}
\multiput(75,55)(0,4){4}{\circle*{2}}
\put(75,80){\vector(0,1){47}}
\put(75,80){\vector(0,1){3}}
\put(75,91){\circle*{4}}
\put(75,109){\circle*{4}}
\multiput(75,55)(4,-2){4}{\circle*{2}}
\put(95,45){\vector(2,-1){44}}
\put(95,45){\vector(2,-1){3}}
\put(91,47){\circle*{4}}
\multiput(75,55)(-4,-2){12}{\circle*{2}}
\put(59,47){\circle*{4}}
\put(43,39){\circle*{4}}
\put(27,31){\circle*{4}}
\put(23,29){\vector(-2,-1){12}}
\put(23,29){\vector(-2,-1){3}}
\put(59,47){\circle*{4}}
\put(71,43){$\scriptscriptstyle{D_4}$}
\put(80,70){$\scriptscriptstyle{S_4}$}
\put(80,88){$\scriptscriptstyle{A_4}$}
\put(80,106){$\scriptscriptstyle{A_4}$}
\put(73,132){$\scriptscriptstyle{3}$}
\put(3,17){$\scriptscriptstyle{2}$}
\put(143,17){$\scriptscriptstyle{4}$}
\put(92,52){$\scriptscriptstyle{D_4}$}
\put(48,52){$\scriptscriptstyle{D_4}$}
\put(32,44){$\scriptscriptstyle{D_2}$}
\put(16,36){$\scriptscriptstyle{D_2}$}
\put(79,57){$\scriptscriptstyle{v'_1}$}
\put(60,39){$\scriptscriptstyle{v'_0}$}
\put(6,32){$\scriptscriptstyle{\varepsilon'_1}$}
\put(10,120){$\scriptstyle{T^{\ast}_{G_1}}$}
\end{picture}
\\
\textrm{Figure 9}
\end{array}
$$
The map $\varrho^0_1$ maps each half-line starting at $v_1$ in
$T^{\ast}_{G_0}$ to the half-line $[v'_1,\varepsilon'_1[$.
Let $G_2$ be isomorphic to $D_{2l}$ with $G_1\bigcap G_2=G_0$ taken as
in \ref{para-5-1}. Then the map $\varrho^0_2$ maps $T^{\ast}_{G_0}$
injectively onto $T^{\ast}_{G_2}$. The center of
$T^{\ast}_{G_2}$ is fixed by $D_2$ if $l$ is odd, or by $D_4$ otherwise.
The center $v_1$ in $T^{\ast}_{G_0}$ is mapped to $v'_1$ fixed 
by $D_4$; if $l$ is even, then the conditon (2) in Lemma \ref{lem-finite}
is not satisfied at $v_1$, since these two $D_4$'s at $v'_1$ and the
center in $T^{\ast}_{G_2}$ are not comparable. Hence $l$ must be 
odd, and in this case, one can verify that the resulting tree of groups
is $\ast$-admissible.
Also, the last argument shows that the first tree in Figure 10 is
$\ast$-admissible, which is simply obtained by two replacements of
mirrors.
$$
\begin{array}{c}
\setlength{\unitlength}{.7pt}
\begin{picture}(150,140)(0,0)
\put(75,55){\circle*{4}}
\put(75,73){\circle*{4}}
\multiput(75,55)(0,4){4}{\circle*{2}}
\put(75,80){\vector(0,1){47}}
\put(75,80){\vector(0,1){3}}
\put(75,91){\circle*{4}}
\put(75,109){\circle*{4}}
\multiput(75,55)(4,-2){4}{\circle*{2}}
\put(95,45){\vector(2,-1){44}}
\put(95,45){\vector(2,-1){3}}
\put(91,47){\circle*{4}}
\multiput(75,55)(-4,-2){12}{\circle*{2}}
\put(59,47){\circle*{4}}
\put(43,39){\circle*{4}}
\put(27,31){\circle*{4}}
\put(23,29){\vector(-2,-1){12}}
\put(23,29){\vector(-2,-1){3}}
\put(59,47){\circle*{4}}
\put(71,43){$\scriptscriptstyle{D_4}$}
\put(80,70){$\scriptscriptstyle{S_4}$}
\put(80,88){$\scriptscriptstyle{A_4}$}
\put(80,106){$\scriptscriptstyle{A_4}$}
\put(73,132){$\scriptscriptstyle{3}$}
\put(0,15){$\scriptscriptstyle{4l+2}$}
\put(143,17){$\scriptscriptstyle{4m}$}
\put(92,52){$\scriptscriptstyle{D_{4m}}$}
\put(48,52){$\scriptscriptstyle{D_4}$}
\put(32,44){$\scriptscriptstyle{D_2}$}
\put(31,24){$\scriptscriptstyle{D_{4l+2}}$}
\end{picture}
\qquad\qquad
\begin{picture}(150,140)(0,0)
\put(75,55){\circle*{4}}
\put(75,73){\circle*{4}}
\multiput(75,56)(0,4){4}{\circle*{2}}
\put(75,80){\vector(0,1){47}}
\put(75,80){\vector(0,1){3}}
\multiput(87,49)(-4,2){4}{\circle*{2}}
\put(95,45){\vector(2,-1){44}}
\put(95,45){\vector(2,-1){3}}
\put(105,40){\circle*{4}}
\put(119,33){\circle*{4}}
\multiput(63,49)(4,2){4}{\circle*{2}}
\put(55,45){\vector(-2,-1){44}}
\put(55,45){\vector(-2,-1){3}}
\put(45,40){\circle*{4}}
\put(31,33){\circle*{4}}
\put(91,47){\circle*{4}}
\put(59,47){\circle*{4}}
\put(78,70){$\scriptscriptstyle{D_{4m}}$}
\put(70,132){$\scriptscriptstyle{4m}$}
\put(3,17){$\scriptscriptstyle{3}$}
\put(143,17){$\scriptscriptstyle{3}$}
\put(92,52){$\scriptscriptstyle{S_4}$}
\put(106,45){$\scriptscriptstyle{A_4}$}
\put(120,38){$\scriptscriptstyle{A_4}$}
\put(48,52){$\scriptscriptstyle{S_4}$}
\put(34,45){$\scriptscriptstyle{A_4}$}
\put(20,38){$\scriptscriptstyle{A_4}$}
\end{picture}
\\
\textrm{Figure 10}
\end{array}
$$
The tree in the right-hand side in Figure 10 with $m=1$ is obtained by
the push-out  
$T^{\ast}_{G_1}\sharp_{T^{\ast}_{G_0}}T^{\ast}_{G'_1}$, where $G_1\cong
S_4$, $G_0\cong D_4$, and $G'_1$ is the twist of $G_1$ by the
non-trivial element in $N(G_0)/(N(G_0)\bigcap N(G_1))\cong Z_2$
(note that $N(D_4)=D_8$). By the suitable replacement of mirrors, we get 
the one with $m\geq 1$.
The center is fixed by $D_4$ if $m$ is odd, or $D_8$ otherwise.

Next consider $G_1\cong A_5$ with a subgroup $G_0\cong D_2$. 
The push-out $T^{\ast}_{G_1}\sharp_{T^{\ast}_{G_0}}T^{\ast}_{G_2}$ with
an embedded group $G_2\cong D_{2l}$ (generated by $G_0$ and an element
of order $2l$ commuting with an involution in $G_0$) yields the
replacement of the mirror of order $2$ by that of order $2l$; but, by a
similar reasoning as above, $l$ must be odd.
Hence we get the tree in the left-hald side of Figure 11. 
$$
\begin{array}{c}
\setlength{\unitlength}{.7pt}
\begin{picture}(150,140)(0,0)
\put(75,55){\circle*{4}}
\put(75,73){\circle*{4}}
\put(75,80){\vector(0,1){47}}
\put(75,80){\vector(0,1){3}}
\multiput(75,55)(4,-2){4}{\circle*{2}}
\put(95,45){\vector(2,-1){44}}
\put(95,45){\vector(2,-1){3}}
\put(75,55){\vector(-2,-1){64}}
\put(75,66){\line(0,-1){11}}
\put(75,66){\vector(0,-1){3}}
\put(91,47){\circle*{4}}
\put(59,47){\circle*{4}}
\put(71,43){$\scriptscriptstyle{A_4}$}
\put(80,70){$\scriptscriptstyle{A_5}$}
\put(73,132){$\scriptscriptstyle{5}$}
\put(3,17){$\scriptscriptstyle{3}$}
\put(142,17){$\scriptscriptstyle{4m+2}$}
\put(91,52){$\scriptscriptstyle{D_{4m+2}}$}
\put(48,52){$\scriptscriptstyle{A_4}$}
\end{picture}
\qquad\qquad
\setlength{\unitlength}{.7pt}
\begin{picture}(150,140)(0,0)
\put(75,55){\circle*{4}}
\put(75,73){\circle*{4}}
\multiput(75,56)(0,4){4}{\circle*{2}}
\put(75,80){\vector(0,1){47}}
\put(75,80){\vector(0,1){3}}
\put(75,55){\line(2,-1){12}}
\put(87,49){\vector(-2,1){3}}
\put(95,45){\vector(2,-1){44}}
\put(95,45){\vector(2,-1){3}}
\put(75,55){\line(-2,-1){12}}
\put(63,49){\vector(2,1){3}}
\put(55,45){\vector(-2,-1){44}}
\put(55,45){\vector(-2,-1){3}}
\put(91,47){\circle*{4}}
\put(59,47){\circle*{4}}
\put(70,43){$\scriptscriptstyle{A_4}$}
\put(79,70){$\scriptscriptstyle{D_{4m+2}}$}
\put(71,132){$\scriptscriptstyle{4m+2}$}
\put(3,17){$\scriptscriptstyle{5}$}
\put(143,17){$\scriptscriptstyle{5}$}
\put(92,52){$\scriptscriptstyle{A_5}$}
\put(48,52){$\scriptscriptstyle{A_5}$}
\end{picture}
\\
\textrm{Figure 11}
\end{array}
$$
The right-hand side of Figure 11 with $m=0$ is obtained from the
tree in the left-hand side by the twist by the non-trivial element in
$N(A_4)/(N(A_4)\bigcap N(A_5))\cong Z_2$; the construction is similar to 
that of the third tree in Figure 8; note that $N(A_4)=S_4$.
Further replacement gives the one with general $m$.

\vspace{1ex}
\paragraph\label{para-2-2}\ 
The following three trees of groups are simply obtained by the
replacement method.
$$
\begin{array}{c}
\setlength{\unitlength}{.7pt}
\begin{picture}(150,140)(0,0)
\put(75,55){\circle*{4}}
\put(75,73){\circle*{4}}
\multiput(75,56)(0,4){4}{\circle*{2}}
\put(75,80){\vector(0,1){47}}
\put(75,80){\vector(0,1){3}}
\put(75,55){\vector(2,-1){64}}
\put(75,55){\vector(-2,-1){64}}
\put(91,47){\circle*{4}}
\put(59,47){\circle*{4}}
\put(71,43){$\scriptscriptstyle{A_4}$}
\put(79,70){$\scriptscriptstyle{D_{2l}}$}
\put(72,132){$\scriptscriptstyle{2l}$}
\put(3,17){$\scriptscriptstyle{3}$}
\put(143,17){$\scriptscriptstyle{3}$}
\put(92,52){$\scriptscriptstyle{A_4}$}
\put(48,52){$\scriptscriptstyle{A_4}$}
\end{picture}
\qquad
\setlength{\unitlength}{.7pt}
\begin{picture}(150,140)(0,0)
\put(75,55){\vector(2,-1){64}}
\put(75,55){\vector(-2,-1){64}}
\put(75,55){\circle*{4}}
\put(75,73){\circle*{4}}
\multiput(75,56)(0,4){4}{\circle*{2}}
\put(75,80){\vector(0,1){47}}
\put(75,80){\vector(0,1){3}}
\put(72,43){$\scriptscriptstyle{2l}$}
\put(80,70){$\scriptscriptstyle{D_m}$}
\put(72,132){$\scriptscriptstyle{m}$}
\put(2,17){$\scriptscriptstyle{2l}$}
\put(142,17){$\scriptscriptstyle{2l}$}
\end{picture}
\qquad
\setlength{\unitlength}{.7pt}
\begin{picture}(150,140)(0,0)
\put(75,55){\circle*{4}}
\put(75,73){\circle*{4}}
\multiput(75,56)(0,4){4}{\circle*{2}}
\put(75,80){\vector(0,1){47}}
\put(75,80){\vector(0,1){3}}
\multiput(75,55)(4,-2){4}{\circle*{2}}
\put(95,45){\vector(2,-1){44}}
\put(95,45){\vector(2,-1){3}}
\multiput(75,55)(-4,-2){4}{\circle*{2}}
\put(55,45){\vector(-2,-1){44}}
\put(55,45){\vector(-2,-1){3}}
\put(91,47){\circle*{4}}
\put(59,47){\circle*{4}}
\put(80,70){$\scriptscriptstyle{D_m}$}
\put(72,132){$\scriptscriptstyle{m}$}
\put(3,17){$\scriptscriptstyle{l}$}
\put(143,17){$\scriptscriptstyle{n}$}
\put(92,52){$\scriptscriptstyle{D_n}$}
\put(48,52){$\scriptscriptstyle{D_l}$}
\end{picture}
\\
\textrm{Figure 12}
\end{array}
$$
Here in the first tree in Figure 12, $l$ must be odd; this follows from
an argument similar to that below Figure 9. 
In the second one also, $m$ is supposed to be odd; if one attempt to
replace it by an even multiple of $m$, then the push-out gives a 
third tree. 
The third one is obtained by replacements of mirrors from
$T^{\ast}_{D_2}$ (hence $l$, $m$, and $n$ are even), and this is
possible if and only if at least two of $l$, $m$, and $n$ are not
multiple of $4$.
The group attached to the certeral vertex is either $D_2$
if all of $l/2$, $m/2$, and $n/2$ is odd, or $D_4$ otherwise.

\vspace{2ex}
\sectioning\label{section-proof}
{\bf Proof of the theorem}

\vspace{1ex}
In this section we prove that the triangle groups obtained in the
previous section are the only possible ones, and complete the proof of
the theorem.
Let $\Gamma$ be a $p$-adic Schwarzian triangle group of Mumford type.
By Corollary \ref{cor-minimal}, the tree $T^{\ast}_{\Gamma}$ is a
tripod.
Let us fix a section $\iota_{\Gamma}\colon T^{\ast}_{\Gamma}\rightarrow
\tree^{\ast}_{\Gamma}$, by which we decorate $T^{\ast}_{\Gamma}$ to be a 
tree of groups.
First we claim:

\vspace{1ex}
\paragraph\label{lem-irreducible}\
{\bf Lemma.}\ {\sl The tree of groups
$(T^{\ast}_{\Gamma},\Gamma_{\bullet})$ is irreducible.}

\vspace{1ex}
\pf
Suffices to invoke Lemma \ref{lem-reducible2}, and the fact that there
is no non-trivial covering over $\P^1_K$ which branches over at most $1$ 
points. 
\qed

\vspace{1ex}
Due to Proposition \ref{pro-directlimit}, the tree of groups
$T^{\ast}_{\Gamma}$ is isomorphic to the direct limit
$\lim_{\longrightarrow}T^{\ast}_{\Gamma_{\bullet}}$.

\vspace{1ex}
\paragraph\label{lem-noncyclic}\
{\bf Lemma.}\ {\sl There exists a vertex $v\in\Vert(T^{\ast}_{\Gamma})$
such that $\Gamma_v$ is not a cyclic group.}

\vspace{1ex}
\pf
If not, all $T^{\ast}_{\Gamma_v}$ are straight lines, and all maps
$T^{\ast}_{\Gamma_v\bigcap\Gamma_w}\rightarrow T^{\ast}_{\Gamma_v}$ are
injective. Hence $\lim_{\longrightarrow}T^{\ast}_{\Gamma_{\bullet}}\cong 
T^{\ast}_{\Gamma}$ must be a straight line.
\qed

\vspace{1ex}
In the sequel, we use the following notation: For
$v_1,v_2\in\Vert(T^{\ast}_{\Gamma})$, we write $G_i=\Gamma_{v_i}$
($i=1,2$) and $G_0=G_1\bigcap G_2$. Let $\varphi$ be the natural
morphism 
$$
\varphi\colon T^{\ast}_{G_1}\sharp_{T^{\ast}_{G_0}}T^{\ast}_{G_2} 
\longrightarrow\lim_{\longrightarrow}T^{\ast}_{\Gamma_{\bullet}}\cong 
T^{\ast}_{\Gamma}
$$
of trees of groups.

\vspace{1ex}
\paragraph\label{para-proof>5}\ 
First we claim that in $p>5$ there is no triangle groups other than
finite subgroups. 
We may assume $G_1$ is not a cyclic group. 

Suppose $G_0$ is not cyclic. Then the map $\varrho^0_1\colon
T^{\ast}_{G_0}\rightarrow T^{\ast}_{G_1}$ has the following properties
(due to \ref{para-generic}):
\begin{condlist}
\item[(1)] $\varrho^0_1$ maps the center of $T^{\ast}_{G_0}$ to
           that of $T^{\ast}_{G_1}$ (since the groups attached to
           outside the center are cyclic).
\item[(2)] The image of $\varrho^0_1$ is a union of half-lines starting
           at the center.
\end{condlist}
Hence, by virtue of Lemma \ref{lem-finite}, either $G_1\subseteq G_2$
or $G_2\subseteq G_1$ must hold (note that in $T^{\ast}_{G_1}$ the
center is fixed by $G_1$).
This implies that the push-out
$T^{\ast}_{G_1}\sharp_{T^{\ast}_{G_0}}T^{\ast}_{G_2}$ equals to either 
$T^{\ast}_{G_1}$ or to $T^{\ast}_{G_2}$.

Therefore, if $\Gamma$ is not finite, then there must be $v_1$ and $v_2$ 
such that $G_0$ is cyclic; by virtue of Lemma \ref{lem-finite}, we may
assume that $G_2$ is not cyclic.
In this case, the push-out
$T^{\ast}_{G_1}\sharp_{T^{\ast}_{G_0}}T^{\ast}_{G_2}$ looks like as
follows (cf.\ Example \ref{exa-pushout}):
$$
\begin{array}{c}
\setlength{\unitlength}{.7pt}
\begin{picture}(150,100)(0,0)
\put(50,50){\vector(-1,1){35}}
\put(50,50){\vector(-1,-1){35}}
\put(50,50){\line(1,0){50}}
\put(100,50){\vector(1,-1){35}}
\put(100,50){\vector(1,1){35}}
\put(50,50){\circle*{4}}
\put(100,50){\circle*{4}}
\put(33,48){$\scriptscriptstyle{G_1}$}
\put(106,48){$\scriptscriptstyle{G_2}$}
\put(70,54){$\scriptscriptstyle{G_0}$}
\end{picture}
\\
\textrm{Figure 13}
\end{array}
$$
If there does not occur a push-out of form
$T\sharp_{T^{\ast}_{G_0}}T^{\ast}_{G_3}$  
in $\lim_{\longrightarrow}T^{\ast}_{\Gamma_{\bullet}}$, then clearly, it
must have more than $3$ ends. But even if it occurs, one can verify that 
the push-out $T\sharp_{T^{\ast}_{G_0}}T^{\ast}_{G_3}$ again has more
than $3$ ends; for example, if $G_3$ is not cyclic, and
$\varrho^0_3\colon T^{\ast}_{G_0}\rightarrow T^{\ast}_{G_3}$ maps the
straight line $T^{\ast}_{G_0}$ to a half-line (necessarily starting at
the center), then the resulting push-out is still like as in
Figure 13 with $G_2$ being replaced by $G_3$.
Other cases are also verified similarly, even easier.

Hence, in any case, it is deduced that the tree $T^{\ast}_{\Gamma}$ must 
have more than $3$ ends, thereby the contradiction.
Therefore, we have proved that in $p>5$ there is no infinite triangle
groups. 

\vspace{1ex}
\paragraph\label{para-proof=5,3,2}\ 
Next we discuss the case $p=5$.
In this case, only $A_5$ is exceptional. In other words, if there is no
$v\in\Vert(T^{\ast}_{\Gamma})$ with $G_v\cong A_5$, then the same
argument as in the previous paragraph shows that $\Gamma$ must be
finite. Hence we may assume that $G_1\cong A_5$.
Again, the argument similar to that in the previous paragraph shows that 
$G_0$ must not be cyclic. Then $G_0$ is isomorphic to either $A_4$,
$D_5$, $D_3$, or $D_2$.
If $G_0$ is not isomorphic to $D_5$, then the image of the map
$\varrho^0_1$ contains a vertex fixed by the whole $G_1$ (cf.\
\ref{para-icosahedral}), and again, $\Gamma$ must be finite. 
Hence $G_0\cong D_5$. The possible isomorphism
class of $G_2$ is $A_5$, $D_{5l}$. 
Now it is an easy combinatorics to check that the push-outs described
in \ref{para-5-1}$\sim$\ref{para-5-3} are the only possible ones.

The case $p=3$ is similar; in this case, $A_4$, $S_4$ and $A_5$ are
exceptional. In each cases, the isomorphism class of $G_0$ is
determined, and the easy observation shows that the trees appeared in
\ref{para-3-1} are the only possible ones.

Also in $p=2$, one can prove that the possible trees are among those in
\ref{para-2-1} and \ref{para-2-2} basically by a similar idea; but the
argument is more involved that the other cases. The possible
non-trivial (isomorphism classes of) pairs $(G_1,G_0)$ are as follows:

\vspace{1ex}
--- $(A_5,A_4)$, $(A_5,D_2)$, 

--- $(S_4,A_4)$, $(S_4,D_4)$, $(S_4,D_2)$,

--- $(A_4,D_2)$, 

--- $(D_{2m+1},Z_2)$, $(D_{2m},D_2)$.

\vspace{1ex}
Here, except for $(D_{2m+1},Z_2)$, cyclic $G_0$ is discarded by the same 
reasoning as in \ref{para-proof>5}; similarly, one sees that the
following cases should be avoided: 
\begin{eqnarray*}
(G_1,G_0,G_2)&=&(A_5,D_2,A_5),(A_5,D_2,S_4),(S_4,D_2,S_4)\ \textrm{(cf.\
Figure 9)},\\
&&(A_4,D_2,A_4),(A_4,D_2,A_5),(A_4,D_2,S_4),(D_{2m+1},Z_2,D_{2n+1}).
\end{eqnarray*}
By a straightforward combinatorics (not very painful but tedious), one
verifies that we have listed all the possible combinations in
\ref{para-2-1} and \ref{para-2-2}.

\vspace{1ex}
\paragraph\label{para-conjugacy}\ 
It remains to prove that the conjugacy class of a $p$-adic Schwarzian
triangle group $\Gamma$ of Mumford type is unique. 
As one finds in the description of $T^{\ast}_{\Gamma}$ in the previous
section, in each $p$, the abstract structure of the tree of groups
$T^{\ast}_{\Gamma}$ is determined by its index $(e_0,e_1,e_{\infty})$; 
in particular, the abstract group structure of $\Gamma$ is determined by 
its index.
Since $T^{\ast}_{\Gamma}$ is a tripod, i.e.\ has exactly three ends, and 
since giving three ends in Bruhat-Tits tree $\tree_K$ determines a
tripod, $\PGL(2,K)$ acts transitively on the set of admissible
embeddings of $T^{\ast}_{\Gamma}$, and hence, on the set of embeddings
of $\Gamma$ in $\PGL(2,K)$.
\qed

\appendix\vspace{2ex}
\sectioning\label{section-appendix}\
{\bf Appendix: Trees of finite groups}

\vspace{1ex}
\paragraph\label{para-appendix}\
This appendix is responsible for detailed description of the tree
$\tree^{\ast}_G$ and the tree of groups $(T^{\ast}_G,G_{\bullet})$
for a finite subgroup $G\subset\PGL(2,K)$.
First let us collect some facts on such subgroups, necessary for later
use, which are either well-known (cf., for example, \cite[Vol.\
II]{Web99}) or are easy to verify:

\vspace{1ex}
(\ref{para-appendix}.1) {\sl Any finite subgroup $G\subset\PGL(2,K)$ is
isomorphic to either a cyclic group (denoted by $Z_m$), a dihedral group
(denoted by $D_m\cong Z_m\rtimes Z_2$), the tetrahedral
group ($\cong A_4$), the octahedral group ($\cong S_4$), or to the
icosahedral group ($\cong A_5$).}

\vspace{1ex}
(\ref{para-appendix}.2) {\sl Two isomoprhic finite subgroups are
conjugate in $\PGL(2,K)$.} 

\vspace{1ex}
(\ref{para-appendix}.3) {\sl Maximal cyclic subgroups of $G$ of same
order comprise a single conjugacy class in $G$ except for the case $G\cong
D_m$ with $m$ even. If
$G=\langle\theta,\varphi\,|\,\theta^m=\varphi^2=(\theta\varphi)^2=1\rangle 
\cong D_m$ with $m$ even, then $\langle\theta\rangle$,
$\langle\varphi\rangle$, and $\langle\theta\varphi\rangle$ give the 
complete system of conjugacy classes of maximal cyclic subgroups.}

\vspace{1ex}
\paragraph\label{para-method}\
{\sl Strategy of description.}\ 
Here is the general strategy for calculating $\tree^{\ast}_G$:

(1) The tree $\tree^{\ast}_G$ is
the minimal one which contains all the mirrors of elements ($\neq 1$) of 
$G$, which are in bijection with maximal cyclic subgroups in $G$.
We therefore first need to know how these mirrors are arranged in
$\tree_K$; the general principle for this will be given in Lemma
\ref{lem-crossratio} below, by which we will see that the necessary data 
are cross-ratios of fixed points. Calculating these values is completely 
elementary, but needs a lot of computation; we will give the complete
list of such data for $G\cong A_4$, $S_4$, and $A_5$ in the next
appendix.
We are thus able to describe the tree $\tree^{\ast}_G$ perfectly.

(2) To describe $(T^{\ast}_G,G_{\bullet})$, we further need to know
the fixed locus in $\tree_K$ of elliptic elements; the fixed locus
contains the mirror, but they do not coincides in general, which has
been already observed in Lemma \ref{lem-fixedlocus}. 

\vspace{1ex}
\paragraph\label{para-crossratio}\
Let $a=(a_0\colon a_1)$, $b=(b_0\colon b_1)$, $c=(c_0\colon c_1)$ and 
$d=(d_0\colon d_1)$ be four pairwise distinct $K$-rational points of
$\P^1_K$. 
We are interested in the arrangement of two apartments $]a,b[$ and
$]c,d[$ in $\tree_K$.
Let us define the cross-ratio
$$
R(a,b;c,d)=
\frac{(a_1c_0-a_0c_1)(b_1d_0-b_0d_1)}{(a_0b_1-a_1b_0)(c_0d_1-c_1d_0)}.
$$

\vspace{1ex}
\paragraph\label{lem-crossratio}\
{\bf Lemma.}\ {\slshape 
Let $\valuation\colon K^{\times}\rightarrow\Z$ be 
the normalized (i.e., $\valuation(\pi)=1$) valuation. 

(1) If $|\valuation(R(a,b;c,d))|=|\valuation(R(b,a;c,d))|=0$, then
$]a,b[$ and $]c,d[$ intersects at exactly one vertex. 

(2) If $|\valuation(R(a,b;c,d))|=|\valuation(R(b,a;c,d))|\neq 0$,
then $]a,b[$ and $]c,d[$ are disjoint with the distance
$|\valuation(R(a,b;c,d))|$. 

(3) If $|\valuation(R(a,b;c,d))|\neq|\valuation(R(b,a;c,d))|$, then
the intersection of $]a,b[$ and $]c,d[$ is the path
$[v(a,b,c),v(b,c,d)]$ of length 
$\max\{|\valuation(R(a,b;c,d))|,|\valuation(R(b,a;c,d))|\}$,
where $v(z_0,z_1,z_2)$ for pairwise distinct three points
$z_0,z_1,z_2\in\P^1(K)$ is the unique vertex lying in the intersection
of all $]z_i,z_j[$ for $i,j=0,1,2$, $i\neq j$.} 

\vspace{1ex}
\pf
First we recall how to calculate $v(z_0,z_1,z_2)$:
Let $Y_i$ ($i=0,1,2$) be a homogeneous coordinate of $z_i$, and choose 
$\alpha_i\in K^{\times}$ such that $\alpha_0Y_0+\alpha_1Y_1+\alpha_2Y_2=0$.
Then $v(z_0,z_1,z_2)$ is the similarity class of
$\O_K\alpha_iY_i+\O_K\alpha_jY_j$ for any $i,j=0,1,2$, $i\neq j$.
By this, it is easily checked that $d(v(a,b,c),v(b,c,d))=
|\valuation(R(a,b;c,d))|$.
Once it is checked, all the statements are clear, since the apartments 
$]a,b[$ and $]c,d[$ either do not intersect or intersect along a path
(see Figure 14). 
\qed

\vspace{-3ex}
\begin{center}
\begin{small}
$$
\setlength{\unitlength}{1pt}
\begin{picture}(350,70)(0,0)
\put(70,20){\vector(1,0){60}}\put(135,18){$d$}
\put(70,20){\vector(-1,0){60}}\put(2,18){$c$}
\put(70,50){\vector(1,0){60}}\put(135,48){$b$}
\put(70,50){\vector(-1,0){60}}\put(2,48){$a$}
\put(70,20){\line(0,1){30}}
\put(70,20){\circle*{4}}\put(70,50){\circle*{4}}
\put(50,5){$v(b,c,d)$}\put(50,60){$v(a,b,c)$}
\put(240,35){\line(1,0){60}}
\put(240,35){\vector(-2,1){30}}\put(240,35){\vector(-2,-1){30}}
\put(300,35){\vector(2,-1){30}}\put(300,35){\vector(2,1){30}}
\put(200,18){$c$}\put(335,18){$d$}
\put(200,48){$a$}\put(335,48){$b$}
\put(240,35){\circle*{4}}\put(300,35){\circle*{4}}
\put(230,20){$v(a,b,c)$}
\put(268,43){$v(b,c,d)$}
\end{picture}
$$
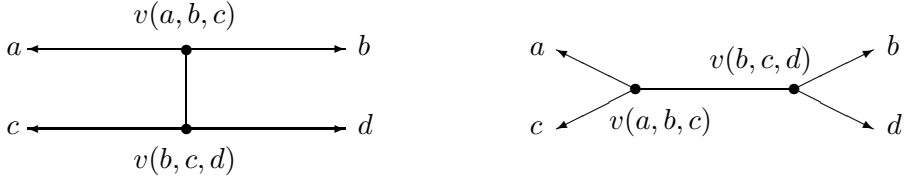
{\sc Figure 14}: Arrangement of the apartments $]a,b[$ and $]c,d[$
\end{small}
\end{center}

\vspace{1ex}
\paragraph\label{para-folding}\
{\sl Folding.}\ 
For describing the quotient tree $T^{\ast}_G$, it is yet more
convenient to have the following notion: Since any element $\gamma\neq 1$
has fixed points in $\Omega_G=\P^{1,\mathrm{an}}_K$, the image of the
mirror $M(\gamma)$ is either a straight-line or a half-line; in
the former case,
$\varrho_G\colon M(\gamma)\rightarrow\varrho_G(M(\gamma))$ is 
an isomorphism, whereas in the latter, it is $2$-to-$1$ except on one
vertex.
If this latter happens, we say that the mirror $M(\gamma)$ is
{\it folded}. 

\vspace{1ex}
(\ref{para-folding}.1)\ {\sl 
The folding of the mirror $M(\gamma)$ occurs if and only if there
exists $\theta\in G$ of order $2$ which interchanges the two fixed
points of $\gamma$, or equivalently, 
$\langle\gamma,\theta\rangle\cong D_m$, where $m$ is the order of
$\gamma$.}

\vspace{1ex}\noindent
Also, it is clear that the folding of mirrors depends only on the
conjugacy classes; hence, in view of (\ref{para-appendix}.3), one can
easily tell which mirror is folded, by only checking existence or
non-existence of dihedral subgroups of $G$.
As a result, we get:

\vspace{1ex}
(\ref{para-folding}.2)\ {\sl If $G\cong D_m$ with $m$ even, $G\cong
S_4$, or $G\cong A_5$, then all mirrors are folded. If $G\cong D_m$ with 
$m$ odd (resp.\ $A_4$), only the mirror of elements of order $m$ (resp.\ 
$2$) is folded. If $G\cong Z_m$, no mirror is folded.
}

\vspace{1ex}
\paragraph\label{para-cyclic}\
{\sl Cyclic case:\ $G\cong Z_n$.}\ 
For any residue characteristic $p$ the tree 
$\tree^{\ast}_G$ consists of only one apartment which is the
mirror $M(\theta)$ of any element $\theta\neq 1$ in $G$.
Since $G$ acts on $\tree^{\ast}_G$ trivially, the quotient tree
$T^{\ast}_G$ also consists of one straight-line whose ends
corresponds to the two points above which $\P^1_K\rightarrow
G\backslash\P^1_K$ ramifies. 

\vspace{1ex}
\paragraph\label{para-convention}\
{\sl Convention.}\ 
In the following paragraphs, we only present the quotient tree
$T^{\ast}_G$ and the stabilizers in pictures for $G$ a non-cyclic
subgroups. 
One can check these by first drawing $\tree^{\ast}_G$ by means of Lemma
\ref{lem-crossratio} and the data in tables in the next appendix;
details are left to the reader (but, as a hint for the careful reader, we
will exhibit in the end of the next appendix the picture of $\tree^{\ast}_G$ 
for $G\cong A_5$ in $p=2$).
The pictures are subject to the following convention: 
\begin{itemize}
\item Solid lines are the images of mirrors, while dotted segments are
      the ones which are not 
images of any mirror (recall that the tree $\tree^{\ast}_G$ is not
in general simply the union of mirrors). \vspace{-1ex}
\item Ends are denoted by the arrow. \vspace{-1ex}
\item If a mirror has the half-line as its
image, then the starting point is denoted by the symbol 
$\begin{picture}
(10,6)(0,0)\put(2,3){\vector(1,0){3}}\put(2,3){\line(1,0){6}}
\end{picture}$, and the half-line starts at the vertex nearest
      it.\vspace{-1ex} 
\item The stabilizers of edges are omitted, since they are simply the
intersection of the stabilizers of their end points. The number $m$
      placed by a vertex or an end indicates that the stabilizer is
      isomorphic to the cyclic group of order $m$.\vspace{-1ex}
\end{itemize}
Also, the unit length $u$ (i.e.\ the distance of neighboring dots) are
given in each picture.

\vspace{1ex}
\paragraph\label{para-generic}\
{\sl Generic cases.}\ 
In the cases $G\cong D_m$ ($m$: odd, $p\neq 2$) and $G\cong A_4$ ($p\neq 
2,3$), the trees $T^{\ast}_G$ look like that in Figure 15; the unique dot 
signifies the vertex fixed by $G$, and the other parts are fixed by
cyclic subgroups of orders in the index in Table 2.
In the cases $G\cong D_m$ ($m$: even, $p\neq 2$), $G\cong S_4$ ($p\neq 
2,3$), and $G\cong A_5$ ($p\neq 2,3,5$), the trees $T^{\ast}_G$ look
like that in Figure 16 with the stabilizers subject to the similar rule.
$$
\setlength{\unitlength}{.7pt}
\begin{picture}(150,150)(0,-10)
\put(75,55){\vector(2,-1){64}}
\put(75,55){\vector(-2,-1){64}}
\put(75,55){\circle*{4}}
\put(75,55){\circle{6}}
\put(75,62){\vector(0,1){65}}
\put(75,62){\vector(0,1){3}}
\put(40,0){{\small {\sc Figure 15}}}
\end{picture}
\qquad\qquad
\begin{picture}(150,150)(0,-10)
\put(81,52){\vector(2,-1){58}}
\put(81,52){\vector(2,-1){3}}
\put(69,52){\vector(-2,-1){58}}
\put(69,52){\vector(-2,-1){3}}
\put(75,55){\circle*{4}}
\put(75,55){\circle{6}}
\put(75,62){\vector(0,1){65}}
\put(75,62){\vector(0,1){3}}
\put(40,0){{\small {\sc Figure 16}}}
\end{picture}
$$

\paragraph\label{para-dihedral}\
{\sl Dihedral case: $G\cong D_m$, $p=2$.}\ 
$$
\setlength{\unitlength}{.7pt}
\begin{picture}(150,150)(0,-10)
\put(75,55){\vector(2,-1){64}}
\put(75,55){\vector(-2,-1){64}}
\put(75,55){\circle*{4}}
\put(75,73){\circle*{4}}
\put(75,73){\circle{6}}
\multiput(75,56)(0,4){4}{\circle*{2}}
\put(75,80){\vector(0,1){47}}
\put(75,80){\vector(0,1){3}}
\put(73,43){$\scriptstyle{2}$}
\put(80,70){$\scriptstyle{D_m}$}
\put(72,132){$\scriptstyle{m}$}
\put(3,17){$\scriptstyle{2}$}
\put(143,17){$\scriptstyle{2}$}
\put(-16,0){{\small{\sc Figure 17}: $m$: odd, $u=e$}}
\end{picture}
\qquad\qquad
\begin{picture}(150,150)(0,-10)
\put(75,55){\circle*{4}}
\put(75,73){\circle*{4}}
\put(75,73){\circle{6}}
\multiput(75,56)(0,4){4}{\circle*{2}}
\put(75,80){\vector(0,1){47}}
\put(75,80){\vector(0,1){3}}
\multiput(75,55)(4,-2){4}{\circle*{2}}
\put(95,45){\vector(2,-1){44}}
\put(95,45){\vector(2,-1){3}}
\multiput(75,55)(-4,-2){4}{\circle*{2}}
\put(55,45){\vector(-2,-1){44}}
\put(55,45){\vector(-2,-1){3}}
\put(91,47){\circle*{4}}
\put(59,47){\circle*{4}}
\put(71,43){$\scriptstyle{D_2}$}
\put(80,70){$\scriptstyle{D_m}$}
\put(72,132){$\scriptstyle{m}$}
\put(3,17){$\scriptstyle{2}$}
\put(143,17){$\scriptstyle{2}$}
\put(92,52){$\scriptstyle{D_2}$}
\put(48,52){$\scriptstyle{D_2}$}
\put(-18,0){{\small{\sc Figure 18}: $m$: even, $u=\frac{e}{2}$}}
\end{picture}
$$

\paragraph\label{para-tetrahedral}\
{\sl Tetrahedral case: $G\cong A_4$.}\ 
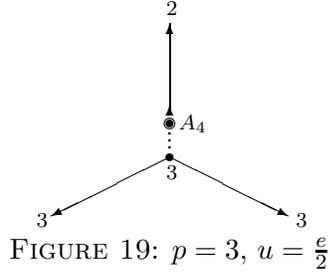
$$
\setlength{\unitlength}{.7pt}
\begin{picture}(150,150)(0,-10)
\put(75,55){\vector(2,-1){64}}
\put(75,55){\vector(-2,-1){64}}
\put(75,55){\circle*{4}}
\put(75,73){\circle*{4}}
\put(75,73){\circle{6}}
\multiput(75,56)(0,4){4}{\circle*{2}}
\put(75,80){\vector(0,1){47}}
\put(75,80){\vector(0,1){3}}
\put(73,43){$\scriptstyle{3}$}
\put(80,70){$\scriptstyle{A_4}$}
\put(73,132){$\scriptstyle{2}$}
\put(3,17){$\scriptstyle{3}$}
\put(143,17){$\scriptstyle{3}$}
\put(-12,0){{\small{\sc Figure 19}: $p=3$, $u=\frac{e}{2}$}}
\end{picture}
\qquad\qquad
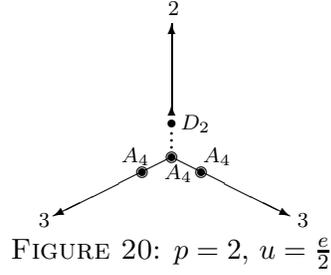
\begin{picture}(150,150)(0,-10)
\put(75,55){\circle*{4}}\put(75,55){\circle{6}}
\put(75,73){\circle*{4}}
\multiput(75,56)(0,4){4}{\circle*{2}}
\put(75,80){\vector(0,1){47}}
\put(75,80){\vector(0,1){3}}
\put(75,55){\vector(2,-1){64}}
\put(75,55){\vector(-2,-1){64}}
\put(91,47){\circle*{4}}\put(91,47){\circle{6}}
\put(59,47){\circle*{4}}\put(59,47){\circle{6}}
\put(71,43){$\scriptstyle{A_4}$}
\put(80,70){$\scriptstyle{D_2}$}
\put(73,132){$\scriptstyle{2}$}
\put(3,17){$\scriptstyle{3}$}
\put(143,17){$\scriptstyle{3}$}
\put(92,52){$\scriptstyle{A_4}$}
\put(48,52){$\scriptstyle{A_4}$}
\put(-12,0){{\small{\sc Figure 20}: $p=2$, $u=\frac{e}{2}$}}
\end{picture}
$$

\paragraph\label{para-octahedral}\
{\sl Octahedral case: $G\cong S_4$.}\ 
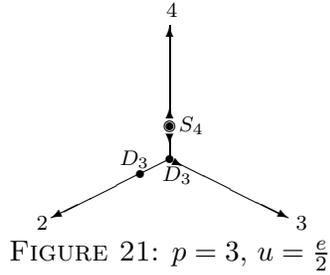
$$
\setlength{\unitlength}{.7pt}
\begin{picture}(150,150)(0,-10)
\put(79,53){\vector(2,-1){60}}
\put(79,53){\vector(2,-1){3}}
\put(75,55){\vector(-2,-1){64}}
\put(75,55){\circle*{4}}
\put(75,73){\circle*{4}}
\put(75,73){\circle{6}}
\put(75,80){\vector(0,1){47}}
\put(75,80){\vector(0,1){3}}
\put(75,66){\line(0,-1){11}}
\put(75,66){\vector(0,-1){3}}
\put(59,47){\circle*{4}}
\put(71,43){$\scriptstyle{D_3}$}
\put(80,70){$\scriptstyle{S_4}$}
\put(73,132){$\scriptstyle{4}$}
\put(3,17){$\scriptstyle{2}$}
\put(143,17){$\scriptstyle{3}$}
\put(48,52){$\scriptstyle{D_3}$}
\put(-12,0){{\small{\sc Figure 21}: $p=3$, $u=\frac{e}{2}$}}
\end{picture}
\qquad\qquad
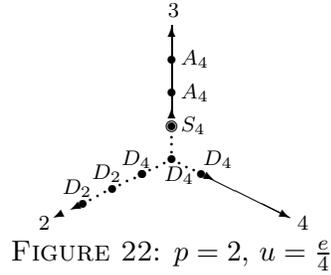
\begin{picture}(150,150)(0,-10)
\put(75,55){\circle*{4}}
\put(75,73){\circle*{4}}
\put(75,73){\circle{6}}
\multiput(75,55)(0,4){4}{\circle*{2}}
\put(75,80){\vector(0,1){47}}
\put(75,80){\vector(0,1){3}}
\put(75,91){\circle*{4}}
\put(75,109){\circle*{4}}
\multiput(75,55)(4,-2){4}{\circle*{2}}
\put(95,45){\vector(2,-1){44}}
\put(95,45){\vector(2,-1){3}}
\put(91,47){\circle*{4}}
\multiput(75,55)(-4,-2){12}{\circle*{2}}
\put(59,47){\circle*{4}}
\put(43,39){\circle*{4}}
\put(27,31){\circle*{4}}
\put(23,29){\vector(-2,-1){12}}
\put(23,29){\vector(-2,-1){3}}
\put(59,47){\circle*{4}}
\put(71,43){$\scriptstyle{D_4}$}
\put(80,70){$\scriptstyle{S_4}$}
\put(80,88){$\scriptstyle{A_4}$}
\put(80,106){$\scriptstyle{A_4}$}
\put(73,132){$\scriptstyle{3}$}
\put(3,17){$\scriptstyle{2}$}
\put(143,17){$\scriptstyle{4}$}
\put(92,52){$\scriptstyle{D_4}$}
\put(48,52){$\scriptstyle{D_4}$}
\put(32,44){$\scriptstyle{D_2}$}
\put(16,36){$\scriptstyle{D_2}$}
\put(-12,0){{\small{\sc Figure 22}: $p=2$, $u=\frac{e}{4}$}}
\end{picture}
$$

\paragraph\label{para-icosahedral}\
{\sl Icosahedral case: $G\cong A_5$.}\ 
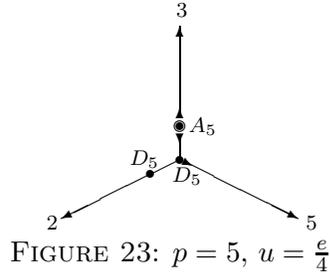
$$
\setlength{\unitlength}{.7pt}
\begin{picture}(150,150)(0,-10)
\put(79,53){\vector(2,-1){60}}
\put(79,53){\vector(2,-1){3}}
\put(75,55){\vector(-2,-1){64}}
\put(75,55){\circle*{4}}
\put(75,73){\circle*{4}}
\put(75,73){\circle{6}}
\put(75,80){\vector(0,1){47}}
\put(75,80){\vector(0,1){3}}
\put(75,66){\line(0,-1){11}}
\put(75,66){\vector(0,-1){3}}
\put(59,47){\circle*{4}}
\put(71,43){$\scriptstyle{D_5}$}
\put(80,70){$\scriptstyle{A_5}$}
\put(73,132){$\scriptstyle{3}$}
\put(3,17){$\scriptstyle{2}$}
\put(143,17){$\scriptstyle{5}$}
\put(48,52){$\scriptstyle{D_5}$}
\put(-17,0){{\small{\sc Figure 23}: $p=5$, $u=\frac{e}{4}$}}
\end{picture}
\qquad\qquad
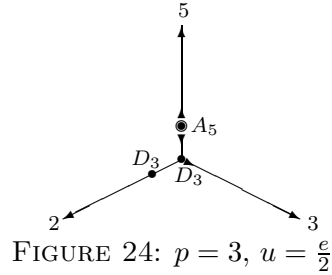
\begin{picture}(150,150)(0,-10)
\put(79,53){\vector(2,-1){60}}
\put(79,53){\vector(2,-1){3}}
\put(75,55){\vector(-2,-1){64}}
\put(75,55){\circle*{4}}
\put(75,73){\circle*{4}}
\put(75,73){\circle{6}}
\put(75,80){\vector(0,1){47}}
\put(75,80){\vector(0,1){3}}
\put(75,66){\line(0,-1){11}}
\put(75,66){\vector(0,-1){3}}
\put(59,47){\circle*{4}}
\put(71,43){$\scriptstyle{D_3}$}
\put(80,70){$\scriptstyle{A_5}$}
\put(73,132){$\scriptstyle{5}$}
\put(3,17){$\scriptstyle{2}$}
\put(143,17){$\scriptstyle{3}$}
\put(48,52){$\scriptstyle{D_3}$}
\put(-17,0){{\small{\sc Figure 24}: $p=3$, $u=\frac{e}{2}$}}
\end{picture}
$$
$$
\setlength{\unitlength}{.7pt}
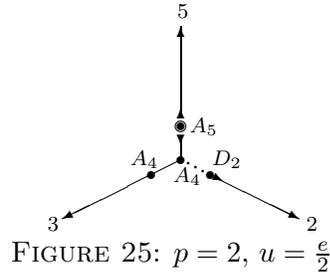
\begin{picture}(150,150)(0,-10)
\put(75,55){\circle*{4}}
\put(75,73){\circle*{4}}
\put(75,73){\circle{6}}
\put(75,80){\vector(0,1){47}}
\put(75,80){\vector(0,1){3}}
\multiput(75,55)(4,-2){4}{\circle*{2}}
\put(95,45){\vector(2,-1){44}}
\put(95,45){\vector(2,-1){3}}
\put(75,55){\vector(-2,-1){64}}
\put(75,66){\line(0,-1){11}}
\put(75,66){\vector(0,-1){3}}
\put(91,47){\circle*{4}}
\put(59,47){\circle*{4}}
\put(71,43){$\scriptstyle{A_4}$}
\put(80,70){$\scriptstyle{A_5}$}
\put(73,132){$\scriptstyle{5}$}
\put(3,17){$\scriptstyle{3}$}
\put(143,17){$\scriptstyle{2}$}
\put(92,52){$\scriptstyle{D_2}$}
\put(48,52){$\scriptstyle{A_4}$}
\put(-17,0){{\small{\sc Figure 25}: $p=2$, $u=\frac{e}{2}$}}
\end{picture}
$$

\sectioning\label{section-appendixB}\
{\bf Appendix: Combinatorial data}

\vspace{1ex}
\paragraph\label{para-appendixB}\
This appendix gives the tables of combinatorial data by which one can
see how the mirrors of elements in finite subgroup in $\PGL(2,K)$ are
arranged in $\tree_K$ so that one can draw the picture of 
$\tree^{\ast}_G$.
The basic principle is as follows: Suppose $G$ is a finite subgroup in
$\PGL(2,K)$ and $\gamma,\theta\in G$ ($\gamma\neq 1,\ \theta\neq 1$).
We assume that fixed points $a,b$ (resp.\ $c,d$) of $\gamma$ (resp.\
$\theta$) are in $\P^1(K)$.
Then the mirror of $\gamma$ (resp.\ $\theta$) is given by $]a,b[$
(resp.\ $]c,d[$).
As we saw in Lemma \ref{lem-crossratio} the correlation between 
$M(g)$ and
$M(h)$ can be calculated by the cross-ratios $R(a,b;c,d)$ and
$R(b,a;c,d)$. 
If we have the complete list of these values for every pair of
generators of maximal cyclic subgroups in $G$, we therefore can perfectly
describe $\tree^{\ast}_G$.

\vspace{1ex}
\paragraph\label{para-conventionB}\
{\sl Convention.}\
For the subgroups $G\cong A_4,S_4,A_5$, which are given in the standard
forms as presented below (this is allowed by (\ref{para-appendix}.2)),
we will give a complete list of maximal cyclic 
subgroups by choosing generators. 
To two of these generators, say $\gamma$ and $\theta$, we associate an
expression $P(\gamma,\theta)=P(\theta,\gamma)$, which is either a
number or an unordered pair of numbers. The meaning of
$P(\gamma,\theta)$ is: 

\vspace{1ex}
--- If it is simply a number, then both $|\valuation(R(a,b;c,d))|$ and 
      $|\valuation(R(b,a;c,d))|$ are equal to
      $|\valuation(P(\gamma,\theta))|$.

--- If it is a pair of numbers, say
      $P(\gamma,\theta)=\{s,t\}$, then, as a set,
      $\{|\valuation(s)|,|\valuation(t)|\}$ coincides with
      $\{|\valuation(R(a,b;c,d))|,|\valuation(R(b,a;c,d))|\}$. 

\vspace{1ex}
\paragraph\label{para-tetra&octa}\
{\sl Case $G\cong A_4$ or $S_4$.}\ 
Since $A_4$ is a subgroup of $S_4$, the calculation can be mixed up.
Let $G$ be a subgroup isomorphic to $S_4$. 
By \cite[\S 73]{Web99} we may assume that $G$ is generated by $\theta$ and
$\chi$ with 
$$
\theta=\left[
\begin{array}{cc}
\sqrt{i}&0\\
0&\frac{1}{\sqrt{i}}
\end{array}
\right],\ \quad
\chi=\left[
\begin{array}{rr}
\frac{1-i}{2}&\frac{1-i}{2}\\
-\frac{1+i}{2}&\frac{1+i}{2}
\end{array}
\right],
$$
where $i$ denotes a primitive $4$-th root of unity.
We set $\omega=\chi\theta\chi\theta^2$.
The group $G$ has three cyclic groups of order $4$
generated by each of $\theta,\chi^2\theta^3,\chi\omega\theta^2$, 
four cyclic groups of order $3$ generated by each of 
$\chi,\theta\chi\theta^3,\theta^2\chi\theta^2,\theta^3\chi\theta$, 
and six cyclic groups of order $2$, not coming from those of order $4$, 
generated by each of 
$\omega,\omega\theta^2,\omega\chi,\chi^2\theta,\omega\chi^2,\theta\chi^2$.
Table 3 presents the $P(g,h)$ for any distinct two of these elements.

The group $G$ has a subgroup isomorphic to $A_4$ generated by
$\theta^2$ and $\chi$.
It has four cyclic groups of order $3$, $\langle\chi\rangle$,
$\langle\theta^2\chi\theta^2\rangle$,
$\langle\theta\chi\theta^3\rangle$,
and $\langle\theta^3\chi\theta\rangle$, and 
three cyclic groups of order $2$, $\langle\theta^2\rangle$,
$\langle\chi\theta^2\chi^2\rangle$,
and $\langle\chi^2\theta^2\chi\rangle$.

\begin{center}
\begin{minipage}{3.7in}
\begin{small}
{\sc Table 3:}\ $\cross(g,h)$ for $g$ and $h$
generators of cyclic subgroups in $G\cong S_4$. Here
$\twothree=\{3,6\}$ and $\threefour=\{4,12\}$.
\end{small}
\end{minipage}

\vspace{2ex}
\footnotesize
\begin{tabular}{|c|cc|cccc|cccccc|}
\hline
&$\chi^2\theta^3$&$\chi\omega\theta^2$&$\chi$&$\theta\chi\theta^3$&
$\theta^2\chi\theta^2$&$\theta^3\chi\theta$&$\omega$&$\omega\theta^2$&
$\omega\chi$&$\chi^2\theta$&$\omega\chi^2$&$\theta\chi^2$\\
\hline
$\theta$&$2$&$2$&$\msq23$&$\msq23$&$\msq23$&
$\msq23$&$2$&$2$&$2\sqrt{2}$&$2\sqrt{2}$&$2\sqrt{2}$&$2\sqrt{2}$\\
$\chi\theta\chi^2$&&$2$&$\msq23$&$\msq23$&
$\msq23$&$\msq23$&$2\sqrt{2}$&$2\sqrt{2}$&$2$&$2$&
$2\sqrt{2}$&$2\sqrt{2}$\\
$\chi^2\theta\chi$&&&$\msq23$&$\msq23$&
$\msq23$&$\msq23$&$2\sqrt{2}$&$2\sqrt{2}$&
$2\sqrt{2}$&$2\sqrt{2}$&$2$&$2$\\
\hline
$\chi$&&&&$\twothree$&$\twothree$&$\twothree$&
$2$&$\m2sq3$&$2$&$\m2sq3$&$2$&$\m2sq3$\\
$\theta\chi\theta^3$&&&&&$\twothree$&$\twothree$&
$\m2sq3$&$2$&$\m2sq3$&$2$&$2$&$\m2sq3$\\
$\theta^2\chi\theta^2$&&&&&&$\twothree$&
$2$&$\m2sq3$&$\m2sq3$&$2$&$\m2sq3$&$2$\\
$\theta^3\chi\theta$&&&&&&&$\m2sq3$&$2$&$2$&$\m2sq3$&
$\m2sq3$&$2$\\
\hline
$\omega$&&&&&&&&$2$&$\threefour$&$\threefour$&
$\threefour$&$\threefour$\\
$\omega\theta^2$&&&&&&&&&$\threefour$&$\threefour$&
$\threefour$&$\threefour$\\
$\omega\chi$&&&&&&&&&&$2$&$\threefour$&$\threefour$\\
$\chi^2\theta$&&&&&&&&&&&$\threefour$&$\threefour$\\
$\omega\chi^2$&&&&&&&&&&&&$2$\\
\hline
\end{tabular}
\end{center}\normalsize

\vspace{1ex}
\paragraph\label{para-icosa}\
{\sl Case $G\cong A_5$.}\ 
By [loc.\ cit.\ \S 74] $G$ is generated by $\theta$ and
$\chi$ with 
$$
\theta=\left[
\begin{array}{cc}
\zeta&0\\
0&1
\end{array}
\right],\ \quad
\chi=\left[
\begin{array}{cc}
\zeta+\zeta^4&1\\
1&-(\zeta+\zeta^4)
\end{array}
\right],
$$
where $\zeta$ is a primitive $5$-th root of unity.
We set $\varphi=\chi\theta^{-1}\chi\theta\chi\theta^{-1}$, which sends
$z\mapsto -1/z$.
The group $G$ has

\vspace{1ex}
--- six cyclic subgroups of order $5$: 
$\langle\theta\rangle$, $\langle\theta\chi\rangle$,
$\langle\theta\chi\theta^3\rangle$,
$\langle\theta^2\chi\theta^2\rangle$,
$\langle\theta^3\chi\theta\rangle$, $\langle\theta^4\chi\rangle$,

--- ten cyclic subgroups of order $3$:
$\langle\theta\chi\theta\rangle$, $\langle\theta^2\chi\rangle$,
$\langle\theta^3\chi\rangle$, $\langle\theta^2\chi\theta\rangle$,
$\langle\theta^3\chi\theta^4\rangle$, 
$\langle\theta\chi\varphi\rangle$, 
$\langle\theta^4\chi\varphi\rangle$,
$\langle\theta^2\chi\varphi\theta^4\rangle$, 
$\langle\theta^2\chi\varphi\theta^2\rangle$,
$\langle\theta^3\chi\varphi\theta\rangle$,

--- and fifteen cyclic subgroups of order $2$: 
$\langle\varphi\theta^i\rangle$,
$\langle\theta^{-i}\chi\theta^i\rangle$,
$\langle\theta^{-i}\chi\varphi\theta^i\rangle$
($i=0,1,2,3,4$).

\vspace{1ex} 
Under these notation we have:

\vspace{1ex}
--- $P(g,h)=\sqrt{5}$ if $g$ and $h$ are generators of
      cyclic groups of order $5$ ($g\neq h$). 

--- $P(g,h)=\sqrt{3}\cdot 5^{1/4}$ if $g$ (resp.\ $h$) is 
      one of the generators of cyclic groups of order $5$ (resp.\
      $3$).

\vspace{1ex}
All the rest are shown in the tables below:

\begin{center}
\begin{minipage}{3.5in}
\begin{small}
{\sc Table 4:}\ $\cross(g,h)$ for $g$ and $h$
generators of cyclic subgroups in $G\cong A_5$ of order $3$. Here
$\twothree=\{3,6\}$.  
\end{small}
\end{minipage}

\vspace{2ex}
\footnotesize
\begin{tabular}{|c|ccccccccc|}
\hline
&$\theta^2\chi$&$\theta^3\chi$&$\theta^2\chi\theta$&
$\theta^3\chi\theta^4$&$\theta\chi\varphi$&$\theta^4\chi\varphi$&
$\theta^2\chi\varphi\theta^4$&$\theta^2\chi\varphi\theta^2$&
$\theta^3\chi\varphi\theta$\\
\hline
$\theta\chi\theta$&
$\twothree$&$\twothree$&$3$&$3$&$\twothree$&$\twothree$&$\twothree$&
$3$&$\twothree$\\
$\theta^2\chi$&&
$3$&$3$&$\twothree$&$\twothree$&$\twothree$&$\twothree$&$\twothree$&$3$\\
$\theta^3\chi$&&&
$\twothree$&$3$&$\twothree$&$\twothree$&$3$&$\twothree$&$\twothree$\\
$\theta^2\chi\theta$&&&&
$\twothree$&$3$&$\twothree$&$\twothree$&$\twothree$&$\twothree$\\
$\theta^3\chi\theta^4$&&&&&
$\twothree$&$3$&$\twothree$&$\twothree$&$\twothree$\\
$\theta\chi\varphi$&&&&&&
$3$&$3$&$\twothree$&$\twothree$\\
$\theta^4\chi\varphi$&&&&&&&
$\twothree$&$\twothree$&$3$\\
$\theta^2\chi\varphi\theta^4$&&&&&&&&
$3$&$\twothree$\\
$\theta^2\chi\varphi\theta^2$&&&&&&&&&
$3$\\
\hline
\end{tabular}\normalsize
\end{center}

\def\a{$\ast$}
\vspace{2ex}
\begin{center}
\begin{minipage}{3.5in}
\begin{small}
{\sc Table 5:}\ $\cross(g,h)$ for $g$ (resp.\ $h$) a
generator of cyclic group subgroups in $G\cong A_5$ of order $5$ (resp.\
$2$). Here \a$=2\cdot 5^{1/4}$.
\end{small}
\end{minipage}

\vspace{2ex}
\footnotesize
\begin{tabular}{|c|ccccccccc}
\hline
&$\varphi$&
$\varphi\theta$&$\varphi\theta^2$&$\varphi\theta^3$&
$\varphi\theta^4$&$\chi$&$\theta\chi\theta^4$&$\theta^2\chi\theta^3$&
$\theta^3\chi\theta^2$\\
\hline
$\theta$&
$2$&$2$&$2$&$2$&$2$&\a&\a&\a&\a\\
$\theta\chi$&
\a&\a&\a&\a&$2$&\a&\a&$2$&\a\\
$\theta\chi\theta^3$&
\a&\a&$2$&\a&\a&$2$&\a&\a&$2$\\
$\theta^2\chi\theta^2$&
$2$&\a&\a&\a&\a&\a&$2$&\a&\a\\
$\theta^3\chi\theta$&
\a&\a&\a&$2$&\a&$2$&\a&$2$&\a\\
$\theta^4\chi$&
\a&$2$&\a&\a&\a&\a&$2$&\a&$2$\\
\hline
\end{tabular}
\end{center}

\begin{center}\footnotesize
\begin{tabular}{cccccc|c|}
\hline
$\theta^4\chi\theta$&$\chi\varphi$&
$\theta\chi\varphi\theta^4$&$\theta^2\chi\varphi\theta^3$&
$\theta^3\chi\varphi\theta^2$&$\theta^4\chi\varphi\theta$&\\
\hline
\a&\a&\a&\a&\a&\a&$\theta$\\
$2$&$2$&$2$&\a&\a&\a&$\theta\chi$\\
\a&\a&$2$&$2$&\a&\a&$\theta\chi\theta^3$\\
$2$&\a&\a&$2$&$2$&\a&$\theta^2\chi\theta^2$\\
\a&\a&\a&\a&$2$&$2$&$\theta^3\chi\theta$\\
\a&$2$&\a&\a&\a&$2$&$\theta^4\chi$\\
\hline
\end{tabular}\normalsize
\end{center}

\vspace{2ex}
\begin{center}
\begin{minipage}{3.5in}
\begin{small}
{\sc Table 6:}\ $\cross(g,h)$ for $g$ (resp.\ $h$) a
generator of cyclic group subgroups in $G\cong A_5$ of order $3$ (resp.\ $2$).
\end{small}
\end{minipage}

\vspace{2ex}
\footnotesize
\begin{tabular}{|c|ccccccccc}
\hline
&$\varphi$&
$\varphi\theta$&$\varphi\theta^2$&$\varphi\theta^3$&
$\varphi\theta^4$&$\chi$&$\theta\chi\theta^4$&$\theta^2\chi\theta^3$&
$\theta^3\chi\theta^2$\\
\hline
$\theta\chi\theta$&
$2$&$\sqrt{6}$&$2\sqrt{3}$&$2\sqrt{3}$&$\sqrt{6}$&$2\sqrt{3}$&$2\sqrt{3}$&
$\sqrt{6}$&$\sqrt{6}$\\
$\theta^2\chi$&
$2\sqrt{3}$&$2\sqrt{3}$&$\sqrt{6}$&$2$&$\sqrt{6}$&$2\sqrt{3}$&$\sqrt{6}$&
$\sqrt{6}$\\
$\theta^3\chi$&
$2\sqrt{3}$&$\sqrt{6}$&$2$&$\sqrt{6}$&$2\sqrt{3}$&$2\sqrt{3}$&$\sqrt{6}$&
$\sqrt{6}$&$2\sqrt{3}$\\
$\theta^2\chi\theta$&
$\sqrt{6}$&$2\sqrt{3}$&$2\sqrt{3}$&$\sqrt{6}$&$2$&$\sqrt{6}$&$\sqrt{6}$&
$2\sqrt{3}$&$2\sqrt{3}$\\
$\theta^3\chi\theta^4$&
$\sqrt{6}$&$2$&$\sqrt{6}$&$2\sqrt{3}$&$2\sqrt{3}$&$\sqrt{6}$&
$2\sqrt{3}$&$2\sqrt{3}$&$2\sqrt{3}$\\
$\theta\chi\varphi$&
$2\sqrt{3}$&$\sqrt{6}$&$\sqrt{6}$&$2\sqrt{3}$&$2$&$2$&$2$&$\sqrt{6}$&
$2\sqrt{3}$\\
$\theta^4\chi\varphi$&
$2\sqrt{3}$&$2$&$2\sqrt{3}$&$\sqrt{6}$&$\sqrt{6}$&$2$&$\sqrt{6}$&
$2\sqrt{3}$&$\sqrt{6}$\\
$\theta^2\chi\varphi\theta^4$&
$\sqrt{6}$&$2\sqrt{3}$&$2$&$2\sqrt{3}$&$\sqrt{6}$&$\sqrt{6}$&$2$&$2$&
$\sqrt{6}$\\
$\theta^2\chi\varphi\theta^2$&
$2$&$2\sqrt{3}$&$\sqrt{6}$&$\sqrt{6}$&$2\sqrt{3}$&$2\sqrt{3}$&
$\sqrt{6}$&$2$&$2$\\
$\theta^3\chi\varphi\theta$&
$\sqrt{6}$&$\sqrt{6}$&$2\sqrt{3}$&$2$&$2\sqrt{3}$&$\sqrt{6}$&
$2\sqrt{3}$&$\sqrt{6}$&$2$\\
\hline
\end{tabular}
\end{center}

\begin{center}\footnotesize
\begin{tabular}{cccccc|c|}
\hline
$\theta^4\chi\theta$&$\chi\varphi$&
$\theta\chi\varphi\theta^4$&$\theta^2\chi\varphi\theta^3$&
$\theta^3\chi\varphi\theta^2$&$\theta^4\chi\varphi\theta$&\\
\hline
$2\sqrt{3}$&$2\sqrt{3}$&$2$&$\sqrt{6}$&$\sqrt{6}$&$2$&$\theta\chi\theta$\\
$\sqrt{6}$&$2$&$\sqrt{6}$&$2$&$\sqrt{6}$&$\sqrt{6}$&$\theta^2\chi$\\
$2\sqrt{3}$&$2$&$\sqrt{6}$&$\sqrt{6}$&$2$&$2\sqrt{3}$&$\theta^3\chi$\\
$2\sqrt{3}$&$\sqrt{6}$&$\sqrt{6}$&$2$&$2\sqrt{3}$&$2$&$\theta^2\chi\theta$\\
$\sqrt{6}$&$\sqrt{6}$&$2$&$2\sqrt{3}$&$2$&$\sqrt{6}$&$\theta^3\chi\theta^4$\\
$\sqrt{6}$&$2\sqrt{3}$&$2\sqrt{3}$&$\sqrt{6}$&$2\sqrt{3}$&$\sqrt{6}$&
$\theta\chi\varphi$\\
$2$&$2\sqrt{3}$&$\sqrt{6}$&$2\sqrt{3}$&$\sqrt{6}$&$2\sqrt{3}$&
$\theta^4\chi\varphi$\\
$2\sqrt{3}$&$\sqrt{6}$&$2\sqrt{3}$&$2\sqrt{3}$&$\sqrt{6}$&$2\sqrt{3}$&
$\theta^2\chi\varphi\theta^4$\\
$\sqrt{6}$&$2\sqrt{3}$&$\sqrt{6}$&$2\sqrt{3}$&$2\sqrt{3}$&$\sqrt{6}$&
$\theta^2\chi\varphi\theta^2$\\
$2$&$\sqrt{6}$&$2\sqrt{3}$&$\sqrt{6}$&$2\sqrt{3}$&$2\sqrt{3}$&
$\theta^3\chi\varphi\theta$\\
\hline
\end{tabular}\normalsize
\end{center}

\def\fourfive{\bigstar}
\vspace{2ex}
\begin{center}
\begin{minipage}{3.5in}
\begin{small}
{\sc Table 7:}\ $\cross(g,h)$ for $g$ and $h$
generators of cyclic group subgroups in $G\cong A_5$ of order $2$.
Here $\fourfive=\{4,4\sqrt{5}\}$ and $\threefour=\{4,12\}$.
\end{small}
\end{minipage}

\vspace{2ex}
\footnotesize
\begin{tabular}{|c|cccccccc}
\hline
&$\varphi\theta$&$\varphi\theta^2$&$\varphi\theta^3$&
$\varphi\theta^4$&$\chi$&$\theta\chi\theta^4$&$\theta^2\chi\theta^3$&
$\theta^3\chi\theta^2$\\
\hline
$\varphi$&
$\fourfive$&$\fourfive$&$\fourfive$&$\fourfive$&$2$&$\fourfive$&$\threefour$&
$\threefour$\\
$\varphi\theta$&&
$\fourfive$&$\fourfive$&$\fourfive$&$\threefour$&$\fourfive$&
$2$&$\fourfive$\\
$\varphi\theta^2$&&&
$\fourfive$&$\fourfive$&$\fourfive$&$\threefour$&$\threefour$&$\fourfive$\\
$\varphi\theta^3$&&&&
$\fourfive$&$\fourfive$&$2$&$\fourfive$&$\threefour$\\
$\varphi\theta^4$&&&&&
$\threefour$&$\threefour$&$\fourfive$&$2$\\
$\chi$&&&&&&
$\threefour$&$\fourfive$&$\fourfive$\\
$\theta\chi\theta^4$&&&&&&&
$\threefour$&$\fourfive$\\
$\theta^2\chi\theta^3$&&&&&&&&
$\threefour$\\
$\theta^3\chi\theta^2$\\
$\theta^4\chi\theta$\\
$\chi\varphi$\\
$\theta\chi\varphi\theta^4$\\
$\theta^2\chi\varphi\theta^3$\\
$\theta^3\chi\varphi\theta^2$\\
\hline
\end{tabular}
\end{center}

\vspace{1ex}
\begin{center}\footnotesize
\begin{tabular}{cccccc|c|}
\hline
$\theta^4\chi\theta$&$\chi\varphi$&
$\theta\chi\varphi\theta^4$&$\theta^2\chi\varphi\theta^3$&
$\theta^3\chi\varphi\theta^2$&$\theta^4\chi\varphi\theta$&\\
\hline
$\fourfive$&$2$&$\threefour$&$\fourfive$&$\fourfive$&$\threefour$&
$\varphi$\\
$\threefour$&$\fourfive$&$\threefour$&$2$&$\threefour$&$\fourfive$&
$\varphi\theta$\\
$2$&$\threefour$&$\fourfive$&$\fourfive$&$\threefour$&$2$&
$\varphi\theta^2$\\
$\threefour$&$\threefour$&$2$&$\threefour$&$\fourfive$&$\fourfive$&
$\varphi\theta^3$\\
$\fourfive$&$\fourfive$&$\fourfive$&$\threefour$&$2$&$\threefour$&
$\varphi\theta^4$\\
$\threefour$&$2$&$\fourfive$&$\fourfive$&$\fourfive$&$\fourfive$&
$\chi$\\
$\fourfive$&$\fourfive$&$2$&$\fourfive$&$\fourfive$&$\fourfive$&
$\theta\chi\theta^4$\\
$\fourfive$&$\fourfive$&$\fourfive$&$2$&$\fourfive$&$\fourfive$&
$\theta^2\chi\theta^3$\\
$\threefour$&$\fourfive$&$\fourfive$&$\fourfive$&$2$&$\fourfive$&
$\theta^3\chi\theta^2$\\
&$\fourfive$&$\fourfive$&$\fourfive$&$\fourfive$&$2$&
$\theta^4\chi\theta$\\
&&$\fourfive$&$\threefour$&$\threefour$&$\fourfive$&
$\chi\varphi$\\
&&&$\fourfive$&$\threefour$&$\threefour$&
$\theta\chi\varphi\theta^4$\\
&&&&$\fourfive$&$\threefour$&
$\theta^2\chi\varphi\theta^3$\\
&&&&&$\fourfive$&
$\theta^3\chi\varphi\theta^2$\\
\hline
\end{tabular}\normalsize
\end{center}

\normalsize
\paragraph\label{para-picture}\ 
{\sl Picture of $\tree^{\ast}_G$ with $G\cong A_5$ and $p=2$.}\
Mirrors are denoted by piecewise linear lines; both ends of them
are indicated by arrows.
Each generator of cyclic subgroups of order $3$ is denoted twice at both
ends of its mirror, while that of cyclic subgroups of order $2$ is
written once beside its mirror.
The mirrors of generators of cyclic subgroups of order $5$ are not
drawn; they simply pass through the center of the picture.
Every unit-segment is of length $e/2$.

\vspace{1ex}
\footnotesize
$$
\setlength{\unitlength}{1pt}
\begin{picture}(340,315)(-170,-150)
\put(0,0){\circle*{4}}
\put(0,0){\line(0,1){120}}
\put(0,0){\line(2,1){108}}
\put(0,0){\line(-2,1){108}}
\put(0,0){\line(2,-3){66}}
\put(0,0){\line(-2,-3){66}}
\put(0,60){\circle*{4}}\put(0,120){\circle*{4}}
\put(54,27){\circle*{4}}\put(108,54){\circle*{4}}
\put(-54,27){\circle*{4}}\put(-108,54){\circle*{4}}
\put(33,-49){\circle*{4}}\put(66,-99){\circle*{4}}
\put(-33,-49){\circle*{4}}\put(-66,-99){\circle*{4}}
\put(0,120){\vector(-1,1){20}}\put(0,120){\vector(1,1){20}}
\put(0,120){\vector(-1,3){10}}\put(0,120){\vector(1,3){10}}
           \put(-44,140){$\theta^3\chi\theta^4$}
           \put(-26,153){$\theta^2\chi$}
           \put(13,153){$\theta\chi\varphi$}
           \put(23,140){$\theta^2\chi\varphi\theta^2$}
\put(108,54){\vector(1,2){14}}\put(108,54){\vector(1,0){30}}
\put(108,54){\vector(1,1){24}}\put(108,54){\vector(3,1){32}}
           \put(113,87){$\theta^3\chi\theta^4$}
           \put(135,80){$\theta^3\chi\varphi\theta$}
           \put(144,65){$\theta^2\chi\varphi\theta^4$}
           \put(142,52){$\theta^2\chi\theta$}
\put(-108,54){\vector(-1,2){14}}\put(-108,54){\vector(-1,0){30}}
\put(-108,54){\vector(-1,1){24}}\put(-108,54){\vector(-3,1){32}}
           \put(-128,87){$\theta\chi\theta$}
           \put(-146,82){$\theta^2\chi$}
           \put(-160,68){$\theta^4\chi\varphi$}
           \put(-168,52){$\theta^2\chi\varphi\theta^4$}
\put(66,-99){\vector(2,-1){28}}\put(66,-99){\vector(0,-1){30}}
\put(66,-99){\vector(1,-1){24}}\put(66,-99){\vector(1,-3){10}}
           \put(98,-117){$\theta\chi\theta$}
           \put(94,-127){$\theta\chi\varphi$}
           \put(76,-138){$\theta^3\chi\varphi\theta$}
           \put(56,-137){$\theta^3\chi$}
\put(-66,-99){\vector(-2,-1){28}}\put(-66,-99){\vector(0,-1){30}}
\put(-66,-99){\vector(-1,-1){24}}\put(-66,-99){\vector(-1,-3){10}}
           \put(-116,-117){$\theta^2\chi\theta$}
           \put(-112,-128){$\theta^4\chi\varphi$}
           \put(-89,-136){$\theta^3\chi$}
           \put(-67,-138){$\theta^2\chi\varphi\theta^2$}
\put(0,60){\line(-2,1){54}}\put(-54,87){\circle*{4}}
\put(0,60){\line(-4,1){56}}\put(-56,74){\circle*{4}}
\put(0,60){\line(-4,3){46}}\put(-46,94){\circle*{4}}
\put(-54,87){\vector(-1,1){24}}\put(-54,87){\vector(-3,1){32}}
\put(-56,74){\vector(-2,1){28}}\put(-56,74){\vector(-1,0){31}}
\put(-46,94){\vector(-4,3){25}}\put(-46,94){\vector(-2,3){17}}
           \put(-94,102){$\varphi\theta^2$}
           \put(-86,117){$\theta^4\chi\theta$}
           \put(-108,78){$\theta^4\chi\varphi\theta$}
\put(54,27){\line(0,1){60}}\put(54,87){\circle*{4}}
\put(54,27){\line(-1,4){14}}\put(40,83){\circle*{4}}
\put(54,27){\line(1,4){14}}\put(68,83){\circle*{4}}
\put(54,87){\vector(-1,4){7}}\put(54,87){\vector(1,4){7}}
\put(40,83){\vector(-1,4){7}}\put(40,83){\vector(-2,3){16}}
\put(68,83){\vector(1,4){7}}\put(68,83){\vector(2,3){16}}
           \put(52,117){$\varphi$}
           \put(25,112){$\chi$}
           \put(78,112){$\chi\varphi$}
\put(-54,27){\line(-2,-3){33}}\put(-87,-23){\circle*{4}}
\put(-54,27){\line(-1,-1){40}}\put(-94,-13){\circle*{4}}
\put(-54,27){\line(-1,-3){18}}\put(-72,-27){\circle*{4}}
\put(-87,-23){\vector(-1,-4){7}}\put(-87,-23){\vector(-1,-1){21}}
\put(-94,-13){\vector(-2,-3){16}}\put(-94,-13){\vector(-3,-2){24}}
\put(-72,-27){\vector(-2,-3){16}}\put(-72,-27){\vector(-1,-4){7}}
           \put(-110,-53){$\varphi\theta^4$}
           \put(-136,-40){$\theta^3\chi\theta^2$}
           \put(-93,-63){$\theta^3\chi\varphi\theta^2$}
\put(-33,-49){\line(2,-3){33}}\put(0,-99){\circle*{4}}
\put(-33,-49){\line(1,-1){40}}\put(7,-89){\circle*{4}}
\put(-33,-49){\line(1,-3){18}}\put(-15,-103){\circle*{4}}
\put(0,-99){\vector(1,-1){21}}\put(0,-99){\vector(1,-2){14}}
\put(7,-89){\vector(1,-1){21}}\put(7,-89){\vector(2,-1){28}}
\put(-15,-103){\vector(2,-3){16}}\put(-15,-103){\vector(1,-3){9}}
           \put(18,-129){$\varphi\theta^3$}
           \put(33,-114){$\theta\chi\theta^4$}
           \put(-13,-140){$\theta\chi\varphi\theta^4$}
\put(33,-49){\line(2,1){54}}\put(87,-22){\circle*{4}}
\put(33,-49){\line(1,1){40}}\put(73,-9){\circle*{4}}
\put(33,-49){\line(4,1){56}}\put(89,-35){\circle*{4}}
\put(87,-22){\vector(1,1){21}}\put(87,-22){\vector(4,1){29}}
\put(73,-9){\vector(3,1){28}}\put(73,-9){\vector(2,3){16}}
\put(89,-35){\vector(3,1){28}}\put(89,-35){\vector(1,0){30}}
           \put(108,-11){$\varphi\theta$}
           \put(94,5){$\theta^2\chi\theta^3$}
           \put(120,-32){$\theta^2\chi\varphi\theta^3$}
\end{picture}
$$
\normalsize
\begin{center}
\begin{small}
{\sc Figure 26}: $\tree^{\ast}_G$ for $G\cong A_5$ and $p=2$
\end{small}
\end{center}

\vspace{1ex}
{\sl Acknowledgments.}\ 
The present work owes much to stimulating conversations with Professor
Yves Andr\'e. He pointed out some missing triangle groups in the
preprint version. 
The author expresses gratitude to him. 
The author thanks Gunther Cornelissen and
Aristeides Kontogeorgis for valuable discussions.
Thanks are also due to Professor Jaap Top who pointed out a mistake in
the first version of the main theorem.
The author thanks Max-Planck-Institut f\"ur Mathematik Bonn for the
nice hospitality.
\normalsize   

\hyphenation{Schottky-kurve}
\begin{small}

{\sc Graduate School of Mathematics, Kyushu University, Hakozaki
Higashi-ku, Fukuoka 812-8581, Japan.}
\end{small}

\begin{thebibliography}{99}
\bibitem[And98]{And98}
              Andr\'e, Y.: {\it $p$-adic orbifolds and $p$-adic triangle 
              groups}, RIMS Kyoto proceedings (S\=urikaisekikenky\=usho
              K\=oky\=uroku) No.\ 
              1073, proceedings of the conference ``rigid geometry and
              group action'' Kyoto, December 1998, 136-159.
\bibitem[Ber90]{Ber90}
              Berkovich, V.\ G.: {\it Spectral theory and analytic
              geometry over non-Archimedean fields}, Mathematical
              Surveys and  Monographs, 33, A.M.S., Providence, 1990.
\bibitem[CKK99]{CKK99}
              Cornelissen, G., Kato, F., Kontogeorgis, A.: {\it Discontinuous
              groups in positive characteristic and automorphisms of Mumford
              curves}, preprint, 1999.
\bibitem[Kat00]{Kat00}
              Kato, F.: {\it Graph theoretic construction of discrete
              groups over $p$-adic fields}, preprint, 2000.
\bibitem[GvP80]{GvP80}
              Gerritzen, L., van der Put, M.: {\it Schottky groups and
              Mumford curves}, Lecture Notes in Math., {\bf 817},
              Springer, Berlin, 1980. 
\bibitem[Her78]{Her78}
              Herrlich, F.: {\it \"{U}ber Automorphismen $p$-adischer
              Schottkykurven}, Dissertation, Bochum, 1978.
\bibitem[Her80a]{Her80a}
              Herrlich, F.: {\it Endlich erzeugbare $p$-adische
              diskontinuierliche Gruppen}, Arch.\ Math. {\bf 35} (1980), 
              505--515.
\bibitem[Her80b]{Her80b}
              Herrlich, F.: {\it Die Ordnung der Automorphismengruppe
              einer $p$-adischen Schottkykurve}, Math.\ Ann.\ {\bf 246}
              (1980), 125--130.
\bibitem[Mag74]{Mag74}
              Magnus, W.: {\it Noneuclidean tesselations and their
              groups}, Acad.\ Press, New York, London, 1974.
\bibitem[Mum72]{Mum72}
              Mumford, A.: {\it An analytic construction of degenerating
              curves over complete local rings}, Compositio Math. {\bf
              24} (1972), 129--174.
\bibitem[Ser80]{Ser80}
              Serre, J-P.: {\it Trees}, Springer-Verlag, Berlin,
              Heidelberg, New York, 1980.
\bibitem[Tak77]{Tak77}
              Takeuchi, K.: {\it Commensurability classes of arithmetic
              triangle groups}, J.\ Fac.\ Sc.\ Univ.\ Tokyo {\bf 24}
              (1977), 201--212. 
\bibitem[V-M80]{V-M80}
              Valentini, R.C., Madan, M.L.: {\it A Hauptsatz of L.E.\
              Dickson and Artin-Schreier extensions}, J.\ Reine Angew.\
              Math.\ {\bf 318} (1980), 156--177.
\bibitem[vdP83]{vdP83}
              van der Put, M.: {\it Etale coverings of a Mumford curve}, 
              Ann.\ Inst.\ Fourier, Grenoble {\bf 33,\ 1} (1983),
              29--52.
\bibitem[vdP97]{vdP97}
              van der Put, M.: {\it The structure of $\Omega$ and its
              quotients $\Gamma\backslash\Omega$}, in Proceedings of the 
              Workshop on ``Drinfeld Modules, Modular Schemes and
              Applications'' (Gekeler, E.-U., van der Put, M., Reversat, 
              M., Van Geel, J.\ ed.), World Scientific, Singapore, New
              Jersey, London, Hong Kong, 1997, 103--112.
\bibitem[Web99]{Web99}
              Weber, H.: {\it Lehrbuch der Algebra}, Chelsea Publishing
              Company, New York, 1899.
\bibitem[Yos87]{Yos87}
              Yoshida, M.: {\it Fuchsian differential equations},
              Aspects of Mathematics Vol.\ E11, Friedr.\ Vieweg \& Sohn, 
              Braunschweig/Wiesbaden, 1987.
\end{thebibliography}
\end{document}